\documentclass[12pt]{amsart}
\usepackage{amssymb}
\usepackage{verbatim}
\usepackage{xspace}
\usepackage{pstricks,pst-node}

\theoremstyle{plain}
\newtheorem{theorem}{Theorem}[section]
\newtheorem{corollary}[theorem]{Corollary}
\newtheorem{proposition}[theorem]{Proposition}
\newtheorem{lemma}[theorem]{Lemma}

\theoremstyle{definition}
\newtheorem{definition}[theorem]{Definition}

\theoremstyle{remark}
\newtheorem{remark}[theorem]{Remark}
\newtheorem*{notation}{Notation}
\newtheorem{example}[theorem]{Example}

\newcommand{\ds}{\displaystyle}

\newcommand{\cstar}{\ensuremath{\text{C}^{*}}-}
\renewcommand{\star}{\ensuremath{{}^{*}}\nobreakdash-\hspace{0 pt}}

\newcommand{\ssadj}[1]{S_{#1}^{\vphantom{*}}S_{#1}^*}
\newcommand{\sadjs}[1]{S_{#1}^*S_{#1}^{\vphantom{*}}}
\newcommand{\SSadj}[2]{S_{#1}^{\vphantom{*}} S_{#2}^*}
\newcommand{\SadjS}[2]{S_{#1}^*S_{#2}^{\vphantom{*}}}

\DeclareMathOperator{\dom}{dom}
\DeclareMathOperator{\ran}{ran}
\DeclareMathOperator{\supp}{supp}

\DeclareMathOperator{\Alg}{Alg}

\DeclareMathOperator{\UHF}{UHF}

\newcommand{\bbA}{\mathbb A}  
\newcommand{\bbC}{\mathbb C}
\newcommand{\bbD}{\mathbb D}

\newcommand{\bbN}{\mathbb N}

\newcommand{\bbR}{\mathbb R}
\newcommand{\bbT}{\mathbb T}
\newcommand{\bbZ}{\mathbb Z}

     \newcommand{\sA}{\mathcal A}
     \newcommand{\sB}{\mathcal B}
     \newcommand{\sC}{\mathcal C}
     \newcommand{\sD}{\mathcal D}
     
     \newcommand{\sF}{\mathcal F}
     \newcommand{\sG}{\mathcal G}

     \newcommand{\sN}{\mathcal N}
     \newcommand{\sO}{\mathcal O}
     \newcommand{\sP}{\mathcal P}

     \newcommand{\sS}{\mathcal S}
     \newcommand{\sT}{\mathcal T}

\begin{document}

\title[Subalgebras of Graph \cstar algebras]{Subalgebras of Graph
  \cstar algebras}
\dedicatory{Dedicated to the memory of Gert Kj\ae rg\aa rd Pedersen}

\author{Alan Hopenwasser}
\address{Department of Mathematics\\
        University of Alabama\\
        Tuscaloosa, AL 35487}
\email{ahopenwa@euler.math.ua.edu}

\author{Justin R. Peters}
 \address{Department of Mathematics\\
          Iowa State University\\
        Ames, IA 50011}
 \email{peters@iastate.edu}
 
 \author{Stephen C. Power}
 \address{Department of Mathematics and Statistics\\
        Lancaster University\\
        Lancaster, LA1 4YF\\
        U.K.}
 \email{s.power@lancaster.ac.uk}

\date{August 13, 2004}

 \begin{abstract}
We prove a spectral theorem for bimodules in the context of 
graph {\ensuremath{\text{C}^{*}}-}algebras.  A bimodule 
 over a suitable abelian algebra is determined 
by its spectrum (i.e., its groupoid partial order)
iff it is generated by the Cuntz-Krieger
partial isometries which it contains iff it is invariant under the
gauge automorphisms.  We study 1-cocycles on the Cuntz-Krieger
groupoid associated with a graph 
{\ensuremath{\text{C}^{*}}-}algebra, obtaining results on when
integer valued or bounded cocycles on the natural AF subgroupoid
extend.  To a finite graph with a total order, we associate a nest 
subalgebra of the graph {\ensuremath{\text{C}^{*}}-}algebra and 
then determine its spectrum.  This is used to
investigate properties of the nest subalgebra. We give a
characterization of the partial isometries in a graph 
{\ensuremath{\text{C}^{*}}-}algebra which normalize a natural
diagonal subalgebra and use this to show that gauge invariant
generating triangular subalgebras are classified by their spectra.
 \end{abstract}
\maketitle

\section{Introduction} \label{s:intro}

Groupoid techniques (``coordinatization'') play a major role in the
study of non-self-adjoint subalgebras of \cstar algebras.  The primary
focus of this approach has been on subalgebras of AF \cstar algebras.
In this paper we apply groupoid techniques to the study of subalgebras
of another extremely important class of \cstar algebras: the graph
\cstar algebras.  We develop a spectral theorem for bimodules which
differs somewhat from the similar theorem for AF \cstar algebras. 
Cocycles are a vital tool in the study of analytic subalgebras 
of AF \cstar algebras; accordingly,  we
investigate cocycles in the Cuntz-Krieger groupoid context.
We also apply our spectral theorem for bimodules to study nest
subalgebras of graph \cstar algebras.  Classification of triangular
subalgebras by their spectra is a central result in the AF context.
We extend this result to the graph \cstar algebra context via a
characterization of normalizing partial isometries which is similar
to, and depends on, the AF analog.

Graph \cstar algebras are constructed from directed graphs.  We shall
make one minor modification in the usual notation for this process:
when concatenating edges to form paths we will read right to left.  We
do this because edges (and paths) correspond to partial isometries in
the graph \cstar algebra and composition of partial isometries is
always read from right to left.  This forces some changes in
terminology from what appears elsewhere in the graph \cstar algebra
literature (relevant changes are mentioned in section~\ref{s:prelim}),
but we believe that it is worth paying this 
small price to make some of the
proofs more natural.  Furthermore, our conventions are in conformity
with the ones now in use in the study of higher rank graph
\cstar algebras and in the study of quiver algebras. 
 Section~\ref{s:prelim} also provides some background material needed
 for the proof of the spectral theorem for bimodules.

Graph \cstar algebras are groupoid \cstar algebras, as shown 
in~\cite{MR98g:46083}.  Since we make substantial use of the groupoid,
and in order to establish terminology, we sketch this construction 
in section~\ref{s:grpoid}.  The bimodules which appear in the spectral
theorem for bimodules are bimodules over a natural abelian subalgebra
of the graph--groupoid \cstar algebra.  From the graph point of view,
this is the \cstar subalgebra generated by all the initial and final
projections of the partial isometries associated with paths (the
Cuntz-Krieger partial isometries).  From the groupoid point of view,
this abelian algebra is the algebra of
continuous functions (vanishing at
infinity) on the space of units.  This abelian
algebra need not be maximal; in section~\ref{s:masa} we show that it
is maximal abelian if, and only if, every loop in the graph has an
entrance. 

In order to define the spectrum 
of a bimodule, we need to be able to view all elements of the groupoid
\cstar algebra as functions on the groupoid.  This is possible for
$r$-discrete groupoids (and all the groupoids in this paper are
$r$-discrete) when they are amenable.  It is proven 
in~\cite{MR98g:46083} and ~\cite{MR1962477} that path space groupoids
are always amenable, so spectral techniques are readily available to us.

The spectral theorem for bimodules was first proven by Muhly and
Solel~\cite{MR90m:46098} for groupoids which are $r$-discrete and
principal.  The groupoids which arise from graph \cstar algebras are
$r$-discrete but, in general, not principal.  In the $r$-discrete
principal groupoid context, every bimodule is determined by its
spectrum.  This is false for graph \cstar algebras.  (It is false even
for the Cuntz algebra $O_n$.)  In sections \ref{s:stb1}~and \ref{s:stb2}
we provide two conditions, each of
which is necessasry and sufficient for a bimodule $\sB$ to be
determined by its spectrum.  One condition is that $\sB$ is determined
by the Cuntz-Krieger partial isometries which it contains; the other
is that $\sB$ is invariant under the gauge automorphisms.  As it
happens, the fact that these two conditions are equivalent to each
other can be proven without use of the groupoid model. 
We prove the equivalence of these two conditions in
section~\ref{s:stb1}, which appears before the description of the groupoid
model, and we give the full spectral theorem for bimodules in
section~\ref{s:stb2}. The argument in section~\ref{s:stb1} appeals
only to the spectral bimodule theorem in the AF \cstar algebra 
case. (See~\cite{MR94g:46001}, for example.)
In section~\ref{s:stb3} we extend the
spectral theorem for bimodules by showing that we can replace the
gauge automorphisms by the one parameter automorphism group naturally
associated with any locally constant real valued cocycle satisfying a
mild technical constraint.

Analytic subalgebras play a major role in the study of subalgebras of
AF \cstar algebras.  Analytic subalgebras are most conveniently
described in terms of cocycles on the AF groupoid.  Two special
classes of cocycles of particular importance are the integer valued
cocycles and the bounded cocycles. The Cuntz-Krieger groupoids which
arise from range finite graphs share some, but not all, of the
properties of AF groupoids.  This results in interesting differences
between the cocycle theories in the two contexts.  In
section~\ref{s:cocycles} we introduce techniques for studying cocycles
on the Cuntz-Krieger groupoid and apply these techniques to
investigate both bounded and integer valued cocycles.  Every
Cuntz-Krieger groupoid contains a natural AF subgroupoid;
section~\ref{s:cocycles} is particularly concerned with the question
of when a cocycle on the AF subgroupoid extends to a cocycle on the
whole groupoid. 

The third author (Power) initiated the study of nest subalgebras of
Cuntz \cstar algebras in~\cite{MR86d:47057} in 1985.  This topic then
lay dormant until the first two authors (Hopenwasser and Peters)
revisited the topic using groupoid techniques in~\cite{hop_peters}.
It turns out that everything which was done for nest subalgebras of
Cuntz \cstar algebras can be extended to the graph \cstar algebra
context (for a finite graph), provided that a suitable order is
imposed on the edges of the graph.  Definitions of an ordered graph
and of an associated nest and nest algebra are given in
section~\ref{s:nsa}.  We characterize the Cuntz-Krieger partial
isometries in the nest algebra and, in turn, the spectrum of the nest
algebra.  This enables us to deduce several results about these nest
subalgebras of graph \cstar algebras; for example, the radical is
equal to the closed commutator ideal.

In~\cite{MR91e:46078}, it was shown that the triangular 
subalgebras $\sA$ of AF C*-algebras $\sB$ for which
$\sA\cap  \sA^*$ is a standard AF masa are classified up to 
isometric isomorphism by their associated topological binary relation, 
or spectrum.
This  reduction of the issue of  isomorphism for TAF algebras to
that of classifying their groupoid partial orders has proven to be
a standard tool for the classification of many families.
We shall obtain an analogous reduction  for  triangular
subalgebras $A$ of a wide class of 
graph C*-algebras where $\sA \cap \sA^*$ is the standard 
masa determined by the generators of $C^*(G)$. As in~\cite{MR91e:46078}
the key step for the proof is the identification of the partial 
isometries 
in $C^*(G)$ which normalise $\sD$ as the elements $v$ for which
\[
\|pvq\| = 0 \text{ or } 1, \text{ for all projections }  p,q \in \sD.
\]
We  obtain this characterisation in section~\ref{s:normalizing} and
apply it to  gauge invariant triangular subalgebras 
in section~\ref{s:triangular}. 

Recall that the tensor (or quiver) algebras of directed graphs
correspond to
the norm closed nonselfadjoint
subalgebras of graph \cstar algebras generated by the 
Cuntz-Krieger generators and 
that for various forms of isomorphism these algebras are known to be
in bijective correspondence with their underlying graphs.  
(See~\cite{math.OA/0309394,math.OA/0309420,math.OA/0309363}.)
We remark that the algebras studied here, being bimodules over the
canonical masa, are quite distinct from these algebras and present
more subtle problems of isomorphism type.

\section{Preliminaries} \label{s:prelim}

Let $G=(V,E,r,s)$ be a directed graph.  As usual, $V$ denotes the set
of vertices and $E$ the set of edges.  The range and source maps are
$r$ and $s$. In this paper, we shall modify slightly the usual procedure
(as it appears in most of the literature) for
associating a graph \cstar algebra $C^*(G)$ to $G$.  (As a
consequence, the description of the groupoid underlying $C^*(G)$ will
also be slightly modified).  This minor change is just notational: a
finite
path $\alpha = \alpha_1 \dots \alpha_n$ is a finite sequence of edges,
or a word, 
which satisfies $r(\alpha_{i+1}) = s(\alpha_i)$ for $i = 1, \dots ,
n-1$.  Infinite paths will be infinite sequences with the same
condition for all~$i$.   Edges and finite paths in $G$ correspond to partial 
isometries in the graph \cstar algebra; with this notational change
the path
$\alpha_1 \alpha_2$, for example, 
corrresponds to $S_{\alpha_1} S_{\alpha_2}$.  
 This notational change will result in modification
of some of the usual conditions concerning graphs which appear in the
literature; for example, the condition that every loop has an exit
will be replaced by the condition that every loop has an entrance; no
sinks will be replaced by no sources, etc.  Although we are deviating
from the usual terminology used in most of the literature on graph
\cstar algebras, we are in conformity with the conventions used for
 higher rank graph \cstar algebras (e.g. in 
\cite{MR2001b:46102}) and also for free semigroup(oid)
algebras and quiver algebras.

Throughout this paper we denote the set of finite paths from $G$ by
$F$ and the set of infinite paths by $P$.  
Range and source maps are defined on $F$ as
follows: if $\alpha = \alpha_1 \dots \alpha_n$ then
$r(\alpha) = r(\alpha_1)$ and $s(\alpha) = s(\alpha_n)$.  Due to our
choice of notation for paths, only the range map can be defined on
$P$; this we do in the obvious way.  Also, if
$\alpha = \alpha_1 \dots \alpha_n \in F$ then the \emph{length} of
$\alpha$ (which is $n$) is denoted by $|\alpha|$.

 We assume that the graph $G$ satisfies the property that
$r^{-1}(v)$ is a finite set, for each vertex $v$.
When this property is satisfied, we say that $G$ is
\emph{range finite}.  (This corresponds to `row finite' in the
literature on graph \cstar algebras.)
The graph
\cstar algebra $C^*(G)$ associated with $G$
is the universal \cstar algebras generated by
a set of partial isometries $\{S_e\}_{e\in E}$ which satisfy the
Cuntz-Krieger relations:
\begin{displaymath}
  S_e^* S_e^{\vphantom{*}} = \sum_{\{f \mid r(f) = s(e)\}} 
S_f^{\vphantom{*}} S_f^*.
\end{displaymath}
(This  minor variation on the usual Cuntz-Krieger relations 
is made to conform to our notation for paths.)  Since we assume
throughout this paper that the graph has no sources, we do not need to 
explicitly include a projection for each vertex.  (If $v$ is a vertex, 
there is an edge $e$ with $s(e)=v$ and $P_v = \sadjs{e}$.)

For any finite path $\alpha$, let 
$S_{\alpha} = S_{\alpha_1} \dots S_{\alpha_k}$.  The Cuntz-Krieger
relations imply that any product of the generators  and their
adjoints can be written in the form $\SSadj{\alpha}{\beta}$.
These are the \emph{Cuntz-Krieger partial isometries} in $\sA$

If $\{S_e\}$ are  Cuntz-Krieger generators for $\sA$ and 
if $z$ is a complex number of absolute value one, then $\{zS_e\}$
is another Cuntz-Krieger family which generates $\sA$.  By the
universality of $\sA$, there is an automorphism $\gamma_z$ such
that $\gamma_z(S_e) = zS_e$, for all edges $e$.  These are the gauge
automorphisms of $\sA$.  Note that for any Cuntz-Krieger partial
isometry $\SSadj {\alpha}{\beta}$, we have
$\gamma_z (\SSadj {\alpha}{\beta}) = z^{|\alpha| - |\beta|}
\SSadj {\alpha}{\beta}$.

The gauge automorphisms are used in~\cite{MR2001k:46084} to determine
when the \cstar algebra generated by a representation of the graph $G$
is isomorphic to the graph \cstar algebra.  As part of that analysis  the
authors identify the fixed point algebra of the gauge automorphisms as
the natural AF subalgebra of $\sA$ and describe a faithful
projection of $\sA$ onto the fixed point algebra.  It is clear that
any Cuntz-Krieger partial isometry $\SSadj {\alpha}{\beta}$ with
$|\alpha| = |\beta|$ is in the fixed point algebra of the gauge
automorphisms.  In fact, these partial isometries generate the fixed
point algebra, which we shall denote by $\sF$.  It is proven in
\cite{MR2001k:46084} than $\sF$ is an AF \cstar algebra.

The projection from $\sA$ onto $\sF$ described in
\cite{MR2001k:46084} is the usual expectation:
\begin{displaymath}
 \Phi_0(f) = \int_{\bbT} \gamma_z (f)\,dz.
\end{displaymath}
This is positive, has norm 1, and is faithful in the sense that
$\Phi_0(f^*f) = 0$ implies that $f = 0$. 

Let $B^*(G)$ denote the \star algebra
 generated by  $\{S_e \mid e \in E \}$,
the Cuntz-Krieger generators 
of $\sA$.  So, $B^*(G)$
is just the linear span of the Cuntz-Krieger
partial isometries.  
 If $a \in B^*(G)$, then  $a$ has an
expansion as a finite sum
\begin{displaymath}
  a = \sum_m \sum_{|\lambda|-|\mu|=m} a_{\lambda \mu}
\SSadj {\lambda}{\mu}.
\end{displaymath}
While this expansion is not unique, 
 each term of the form
$\ds \sum_{|\lambda|-|\mu|=m} a_{\lambda \mu}
\SSadj {\lambda}{\mu}$ is completely determined by $a$.   Given
$a$ represented as above, let
\begin{displaymath}
  \Phi_m(a) = \sum_{|\lambda|-|\mu|=m} a_{\lambda \mu}
\SSadj {\lambda}{\mu}.
\end{displaymath}
Since for any $\alpha$ and $\beta$, 
$\gamma_z (\SSadj {\alpha}{\beta}) = z^{|\alpha|-|\beta|}
\SSadj {\alpha}{\beta}$, we have
\begin{displaymath}
  \int_{\bbT} z^{-m} \gamma_z (\SSadj {\alpha}{\beta})\,dz =
\begin{cases}
\SSadj {\alpha}{\beta}, \quad &\text{if } |\alpha|-|\beta| = m, \\
0, &\text{if } |\alpha|-|\beta| \ne m.
\end{cases}
\end{displaymath}
It follows that 
\begin{displaymath}
  \Phi_m (a) = \int_{\bbT} z^{-m} \gamma_z(a)\,dz
\end{displaymath}
for all $a \in B^*(G)$.
Since $\Phi_m$ is well-defined,
linear, and norm decreasing on $B^*(G)$; it 
extends to  all of $\sA$.

Now fix $a \in \sA$ and consider the function
$f \colon \bbT \to \sA$ given by
$f(z) = \gamma_z(a)$.  The Fourier coefficients for $f$ are just the
elements $\Phi_m(a)$ of $\sA$ and we have the Fourier series
$f \sim \sum_{m \in \bbZ} \Phi_m(a) z^m$.  While the infinite sum need
not be convergent, the Cesaro means converge uniformly to $f$.  Since
$f(1) = a$, we obtain the fact that $a$ is in the closed linear span
of the elements $\Phi_m(a)$.  Thus we  have the formal series
\begin{displaymath}
  a \sim \sum_{m \in \bbZ} \Phi_m(a).
\end{displaymath}
with a Cesaro convergence of the series.
We reiterate that $\Phi_0$ maps $\sA$ onto the core AF subalgebra $\sF$.

\section{The Spectral Theorem for Bimodules -- Part I} \label{s:stb1}

A portion of the spectral theorem for bimodules can be proven without
reference to the groupoid model.   The
full theorem will appear in section~\ref{s:stb2} and a further
extension is given in section~\ref{s:stb3}.

Let $\sD$ be the abelian subalgebra of $\sA$ generated by all
projections of the form $\ssadj{\alpha}$ and $\sadjs{\alpha}$.
 This is clearly a subalgebra of the core AF
  algebra $\sF$.  In general, $\sD$ need not be maximal abelian in 
$\sA$ (though it will be maximal abelian in $\sF$).  We discuss when
$\sD$ will be maximal abelian in $\sA$ in section~\ref{s:masa}.

\begin{theorem} \label{t:stb_prelim}
Let $G$ be a range finite graph with no sources.  Let
$\sB \subseteq \sA$ be a bimodule over $\sD$.  Then $\sB$ is generated
by the Cuntz-Krieger partial isometries which it contains if, and only
if, it is invariant under the gauge automorphisms.
\end{theorem}

\begin{proof}
It is trivial that a bimodule generated by its Cuntz-Krieger partial
isometries is invariant under the gauge automorphisms, so we need only
prove the converse.

 Let $\sB$ be a gauge invariant bimodule over $\sD$.  First note
 that for each $m$, $\Phi_m(\sB) \subseteq \sB$.  For each path 
$\nu \in F$, let
\begin{displaymath}
  \sB^{\nu} = \{b \in \sF \mid S_{\nu} b \in \sB \}.
\end{displaymath}
We claim that $\sB^{\nu}$ is a closed bimodule over $\sD$.  It is
trivial to see that $\sB^{\nu}$ is closed and a 
right bimodule.  Since
$\sD$ is generated by projections of the form $\ssadj {\alpha}$, we
can show that $\sB^{\nu}$ is a 
left bimodule by showing that for each 
$b \in \sB$ and each finite path $\alpha$, 
the element $S_{\nu} (\ssadj {\alpha}) b \in \sB$.
Such an element is non-zero when
$\ssadj{\alpha} \leq \sadjs{\nu}$, and in this case
\begin{displaymath}
S_{\nu}^{\vphantom{*}} (\ssadj {\alpha}) b = 
S_{\nu}^{\vphantom{*}} \ssadj{\alpha} \sadjs{\nu}b =
(S_{\nu}^{\vphantom{*}} \ssadj{\alpha} S_{\nu}^{*} ) S_{\nu}^{\vphantom{*}}b.
\end{displaymath}
This is in $\sB$, since 
$S_{\nu}^{\vphantom{*}} \ssadj{\alpha} S_{\nu}^{*} \in \sD$.
Similarly, the spaces
\begin{displaymath}
  \sB_{\nu} = \{ b \in \sF \mid b S_{\nu}^{*} \in \sB \}
\end{displaymath}
are also closed $\sD$-bimodules.  Since $\sD$ is a canonical masa in
the AF algebra $\sF$ and the $\sB^{\nu}$ and $\sB_{\nu}$ are 
$\sD$-bimodules in $\sF$, the spectral theorem for bimodules in the AF
case implies that each of $\sB^{\nu}$ and $\sB_{\nu}$ is spanned by
the matrix unit elements $\SSadj{\alpha}{\beta}$ in 
$\sB^{\nu}$ or $\sB_{\nu}$ (as appropriate) with
$|\alpha| = |\beta|$.  Thus, the spaces $S_{\nu} \sB^{\nu}$ and 
$\sB_{\nu}^{\vphantom{*}}S_{\nu}^*$ are generated by their
Cuntz-Krieger partial isometries.

We claim that it follows that the spaces $\Phi_m(\sB)$ are also
generated by their Cuntz-Krieger partial isometries.
In view of Cesaro convergence and the fact that 
the $\Phi_m(\sB)$ spaces are subspaces of
$\sB$, this implies that $\sB$ is generated by its Cuntz-Krieger
partial isometries.

The claim is elementary to confirm in the case of a finite graph, since
$\Phi_m(\sB)$ is then the finite linear 
span of the spaces $S_{\nu} \sB^{\nu}$
or $\sB_{\nu}^{\vphantom{*}}S_{\nu}^*$
with $|\nu| = m$ or $|\nu| = -m$, as appropriate,
 and the isometries $S_\nu$ have orthogonal ranges.
In general the claim will follow if we show that  $\Phi_m(\sB)$ 
is the closed linear span of these subspaces.

The case when the graph is infinite can be reduced to the finite graph
case as follows.  Recall that the Cuntz-Krieger partial isometries in
$\Phi_m(\sA)$ are precisely the $\SSadj{\alpha}{\beta}$ with
$|\alpha|-|\beta|=m$.  Let $F_n$ be a sequence of finite subsets of
the Cuntz-Krieger partial isometries in $\Phi_m(\sA)$ such that
$\bigcup F_n$ is the set of all Cuntz-Krieger partial isometries
in $\Phi_m(\sA)$. Also, let $P_n$ denote the projection onto the
closed linear span of the ranges of the partial isometries in $F_n$.

 Any element $b$ in $\sA$ can be approximated by a
linear combination of Cuntz-Krieger partial isometries. But
$\Phi_m$ is contractive, acts as the identity on Cuntz-Krieger partial
isometries in $\Phi_m(\sA)$ and maps all other Cuntz-Krieger partial
isometries to $0$; consequently,
when $b \in \Phi_m(\sA)$ it can be approximated  by linear
combinations of 
Cuntz-Krieger partial isometies in $\Phi_m(\sA)$.  In particular, there
is a sequence $a_n \in \Phi_m(\sA)$ such that
$P_na_n = a_n$ and $a_n \to b$.  Now, suppose further that
$b\in \Phi_m(\sB)$.  Since 
$P_nb -b = P_n(b-a_n) + a_n-b$, we have $P_nb \to b$.  Also
$P_nb \in P_n \Phi_m(\sB) = \Phi_m(P_n\sB)$.  By the result for
finite graphs, each $P_nb$ can be approximated by
linear combinations of Cuntz-Krieger
partial isometries in $\Phi_m(P_n\sB)$.  It follows that $b$ can be
approximated by linear combinations of Cuntz-Krieger partial
isometries in $\Phi_m(\sB)$. 
\end{proof}

\section{The Groupoid Model} \label{s:grpoid}

In~\cite{MR98g:46083}, Kumjian, Pask, Raeburn and Renault construct a
locally compact $r$-discrete groupoid $\sG$ such that 
the groupoid \cstar algebra $C^*(\sG)$ is the graph \cstar algebra for
$G$. We sketch below a slightly modified version of this construction.

We shall assume that 
every vertex is the range
of at least one edge.  (The graph has no sources.)
It follows that every edge is part of an
infinite path (notationally, infinite to the right).
 Infinite path space
$P$ is topologized by taking as a basis of open sets the 
following cylinder
sets: for each finite path $\alpha$ of length $k$,
\begin{align*}
Z(\alpha) &= \{ x \in P \mid x_1=\alpha_1, \,\dots,
\,x_k=\alpha_k \} \\
&= \{\alpha y \mid y \in P \text{ and } r(y)=s(\alpha) \}.
\end{align*}
Any two cylinder sets $Z(\alpha)$ and $Z(\beta)$ are either disjoint
or one is a subset of the other.  For example, 
$Z(\alpha) \subseteq Z(\beta)$ precisely when
$\alpha = \beta \alpha'$ for some $\alpha' \in F$ with
$r(\alpha') = s(\beta)$.
The assumption that $G$ is range finite implies that each
$Z(\alpha)$ is a compact set.  The topology on $P$ is then 
locally compact, $\sigma$-compact, totally disconnected and
Hausdorff. It coincides with the relative product topology obtained by
viewing $P$ as a subset of the infinite product of $E$ with
itself.  Path space $P$ will, in due course, be identified with the
space of units for the groupoid $\sG$.

The next step  is to define
an equivalence relation (\emph{shift equivalence}) on $P$.  Shift
equivalence is the union of a sequence of relations, indexed by the
integers.  Let $x,y \in P$ and $k \in \bbZ$.  If there is a
positive integer $N$ such that $x_{i+k} = y_i$ for all $i \geq N$,
then we write $x \sim_k y$. 
 We then say that $x$ and $y$ in $P$ are \emph{shift
  equivalent} if $x \sim_k y$ for some $k \in \bbZ$.

The groupoid is defined to be the set:
\[
\sG = \{(x,k,y) \mid x,y \in P, k \in \bbZ, x \sim_k y\}.
\]
Elements $(x,k,y)$ and $(w,j,z)$ are composable if, and only if, $y=w$;
when this is the  case $(x,k,y)\cdot(y,j,z)=(x,k+j,z)$.  
Inversion is given by
$(x,k,y)^{-1} = (y,-k,x)$.  With these operations $\sG$ is a
groupoid.  The units of $\sG$ all have the form $(x,0,x)$ for 
$x \in P$, allowing the identification of $P$ with the space of
units (which is also denoted by $\sG^0$, as usual).  

There is a natural topology which renders $\sG$ a topological
groupoid.  A basis for this topology can be parameterized by pairs of
finite paths $\alpha$ and $\beta$ which satisfy $s(\alpha)=s(\beta)$.
  For such $\alpha$ and $\beta$, let
\begin{align*}
Z(\alpha,\beta) 
&= \{(x,k,y) \mid x \in Z(\alpha), y \in Z(\beta), k=|\alpha|-|\beta|,
\text{ and } x_{i+k}=y_i \text{ for } i > |\beta| \} \\
&= \{(\alpha z, |\alpha|-|\beta|, \beta z) \mid z \in P,
r(z)=s(\alpha)=s(\beta) \}.
\end{align*}
 
We allow either $\alpha$ or $\beta$ (or both) to be the empty paths.
For example,
\[
Z(\alpha,\emptyset) = \{ (\alpha z, |\alpha|, z) \mid z \in P,
r(z)=s(\alpha) \}.
\]

Two basic open sets, $Z(\alpha,\beta)$ and $Z(\gamma, \delta)$ are
either disjoint or one contains the other.  For example,
$Z(\alpha,\beta) \subseteq Z(\gamma,\delta)$ precisely when there is
$\epsilon \in F$ such that $\alpha = \gamma \epsilon$ and
$\beta = \delta \epsilon$.  The following proposition is essentially
quoted from~\cite{MR98g:46083}.

\begin{proposition} \label{p:groupoid}
 The sets
\[
\{Z(\alpha,\beta) \mid \alpha, \beta \in F, s(\alpha)=s(\beta) \}
\]
form a basis for a locally compact Hausdorff topology on $\sG$.  With
this topology, $\sG$ is a second countable, $r$-discrete locally
compact groupoid in which each $Z(\alpha,\beta)$ 
\/{\rm (}except possibly $Z(\emptyset, \emptyset)$\/{\rm )}
is a compact open
$\sG$-set.  The product topology on the unit space $P$ agrees with
the topology it inherits by viewing it as the subset
$\sG^0 = \{(x,0,x) \mid x \in P \}$ of $\sG$.  The counting
measures form a left Haar system for $\sG$.
\end{proposition}

Kumjian, Pask, Raeburn, and Renault prove that the groupoid \cstar
algebra for $\sG$ is isomorphic to the graph \cstar algebra $C^*(G)$;
this is done by identifying natural Cuntz-Krieger generators
in $C^*(\sG)$ and proving universality.  Recall that the groupoid
\cstar algebra is constructed by providing $C_c(\sG)$, the compactly
supported continuous functions on $\sG$, with a suitable (convolution
style) multiplication, an involution, and a (universal) \cstar norm
and then completing the \star algebra.  In particular, for each edge
$e$, the set 
$Z(e,\emptyset) = \{(ez,1,z) \mid z \in P, r(z) = s(e) \}$
is compact and open;
therefore its characteristic function
$\chi_{Z(e,\emptyset)}$ may be viewed as an element of $C^*(\sG)$.  
Denote this element by $S_e$.  A routine  calculation shows that
the initial space $\sadjs {e}$ is the characteristic function of
$\{(x,0,x) \mid r(x) = s(e)\}$. 
 Another routine calculation shows that for an edge $f$,
$\ssadj {f} = \chi_{Z(f,f)}$.  Now $Z(f,f)$ is a subset of
$\{(x,0,x) \mid r(x) = s(e)\}$ exactly when $r(f) = s(e)$ and, in fact,
$\{(x,0,x) \mid r(x) = s(e)\}$ is the union of all $Z(f,f)$ with 
$r(f) = s(e)$.  Thus, the Cuntz-Krieger relations
\[
\sadjs{e} = \sum_{r(f)=s(e)} \ssadj{f}
\]
hold.  Routine but tedious calculations show that
$\SSadj{\alpha}{\beta} = \chi_{Z(\alpha,\beta)}$.  
The following theorem consists of a combination of parts of
Proposition 4.1 and Theorem 4.2  in~\cite{MR98g:46083}):

\begin{theorem}
Let $G$ be a range finite directed graph with no sources.
With the notation above, $C^*(\sG)$ is generated by 
$\{S_e \mid e \in G\}$ and $C^*(\sG)$ is isomorphic to $C^*(G)$.
\end{theorem}

Throughout the rest of this paper the graph \cstar algebra--groupoid
\cstar algebra determined by the graph $G$ will be denoted by $\sA$.

\section{The Masa Theorem} \label{s:masa}

In the groupoid model
there is a natural abelian subalgebra of $\sA$: the
functions supported on the space of units of $\sG$.  
We shall denote this abelian algebra
by $\sD$. This is, of course, exactly the same abelian
algebra as the one that appears in section~\ref{s:stb1}.
 Based on what happens for $r$-discrete principal
groupoids and for the Cuntz groupoids which model 
the Cuntz algebras $O_n$,
it might be suspected that $\sD$ is always a masa
in $\sA$; however, this is not  the case.  Consider, for
example, the graph which consists of a single vertex and a single edge
$e$.  Then there is only one infinite path, say $x=eee\dots$ and the
unit space consists of the singleton $(x,0,x)$.  The whole groupoid may
be identified with $\bbZ$: $\sG = \{(x,k,x) \mid k \in \bbZ\}$ and
$\sA \cong C(\bbT)$ while $\sD \cong \bbC$.  

For a more interesting example, let $G$ consist of a single loop with
$n$ edges.  So $E = \{e_1, e_2, \dots, e_n \}$ with $r(e_j) =
s(e_{j-1})$ for $j = 2, \dots, n$ and $r(e_1) = s(e_n)$.
In this case the graph \cstar algebra is
$M_n(C(\bbT))$.  The algebra
$\sD$ corresponds to the diagonal
matrices with scalar entries, which is not a masa.

We will prove  that $\sD$ is a masa in $\sA$ if, and
only if every loop has an entrance.
   This condition says that
 if the finite path $\alpha_1 \dots \alpha_n$ satisfies
$r(\alpha_1) = s(\alpha_n)$, then there is an index $j$ and an edge
$\beta$ such that $\beta \neq \alpha_j$ and $r(\beta) = r(\alpha_j)$.
This condition was used earlier in the literature (expressed as
``every loop has an exit,'' of course).
In~\cite{MR99i:46049} and 
in~\cite{MR2001k:46084},
for example, it is shown that when this
condition holds the \cstar algebra
generated by any system of Cuntz-Krieger partial isometries is
isomorphic to the universal graph \cstar algebra.

The isotropy group bundle of $\sG$ is 
$\sG^1 = \{(x,k,y) \mid x=y \}$.  The space of units of $\sG$ is 
$\sG^0 = \{(x,0,x) \mid x \in P\}$
 The following lemma, combined with groupoid amenability and some
 results in~\cite{MR82h:46075}, will yield the masa theorem.

\begin{lemma} \label{l:interior}
Let $G$ be a range finite directed graph with no sources.
  Let $\sG$ be the
associated groupoid.  Then every loop in $G$ has an entrance if, and
only if, $\sG^0$ is the interior of $\sG^1$.
\end{lemma}

\begin{proof}
Assume that every loop has an entrance.  Let $(x,k,x)$ be an element
of $\sG^1$ which is not in $\sG^0$; in other words, assume that 
$k \neq 0$.  We shall show that $(x,k,x)$ can be approximated by
elements of the complement of $\sG^1$.  Since $\sG^0$ is open, this
will show that $\sG^0$ is the interior of $\sG^1$.

Since $k \neq 0$,
there is an
integer $N$ such that for $i \geq N$, $x_{i+k}=x_i$.  Let
$\beta = x_1 \dots x_{N-1}$, a finite path of length $N-1$ and
let $\alpha = x_N \dots x_{N+k-1}$, a finite path of length $k$.  The
condition for shift equivalence assures that
$x_{N+k} \dots x_{N+2k-1} = x_N \dots x_{N+k-1}$, etc.  Thus,
$x = \beta \alpha \alpha \alpha \dots$.  

Since $\alpha$ can be concatenated with itself,  $\alpha$ is a loop.
Write $\alpha = \alpha_1 \dots \alpha_k$, where the
$\alpha_i \in E$.  Since every loop has an
entrance, there is an edge $y_j$ such that $y_j \neq \alpha_j$ and
$r(y_j) = r(\alpha_j)$.  Now let $y = y_j y_{j+1} \dots$ be an
infinite path ending in the edge $y_j$.  The assumption that the graph
has no sources guarantees the existence of such an infinite path.

For each integer $p \geq 1$, let $z^p$ be the infinite path
$\beta \alpha \dots \alpha \alpha_1 \dots \alpha_{j-1} y$, where there
are exactly $p$ copies of $\alpha$.   If $k>0$,
then $z^{p+1} \sim_k z^p$ and if $k<0$ then $z^p \sim_k z^{p+1}$.
Now just observe that
$z^p \neq z^{p+1}$ and that $(z^{p+1},k,z^p)$ or
$(z^p,k,z^{p+1})$, as appropriate, converges to $(x,k,x)$,

For the converse,
assume that $G$ has a loop $\alpha$ with no entrance.  Let $k$ be
the length of this loop and let $x = \alpha\alpha \alpha \dots$.
Then $(x,k,x) \in \sG^1 \setminus \sG^0$ and the singleton set
$\{(x,k,x)\}$ is open in $\sG$.  Thus, $\sG^0$ is not the interior
of $\sG^1$.
\end{proof}

\begin{theorem} \label{t:masa}
Let $G$ be a range finite directed graph with no sources.
Let $\sG$ be the associated groupoid.
  Then $\sD$ is a masa in $\sA$ if, and only if, every loop
has an entrance.
\end{theorem}

\begin{proof}
Results in~\cite{MR98g:46083} and~\cite{MR1962477} establish the
amenability of $\sG$.  (This is proven for locally finite graphs 
in~\cite{MR98g:46083} and extended to range finite graphs -- and
beyond -- in~\cite{MR1962477}.)  It follows that
$C^*(\sG) = C^*_{\text{red}}(\sG)$. Proposition II.4.7
in~\cite{MR82h:46075} says that $\sD$ is a masa in
$ C^*_{\text{red}}(\sG)$ if, and only if, $\sG^0$ is the interior of
$\sG^1$, so this, combined with
 Lemma~\ref{l:interior} yields the theorem.
\end{proof}

\section{The Spectral Theorem for Bimodules -- Part II} \label{s:stb2}

One of the most useful tools in the study of nonselfadjoint
subalgebras of \cstar algebras is the spectral theorem for bimodules of
Muhly and Solel~\cite{MR90m:46098}.  (See also~\cite{MR94i:46075} for a
generalization due to Muhly, Qiu and Solel.)  The theorem as it
appears in these two references is not valid for graph \cstar
algebras; this section is devoted to the proof of a modified version
of the spectral theorem for bimodules which is valid for a broad class
of graph \cstar algebras.  The theorem as it appears here was proven in
the context of the Cuntz algebras $\sO_n$ in~\cite{hop_peters}; the
proof of the general version is similar.  We will give an extension of
the spectral theorem for bimodules in section~\ref{s:stb3}.

Throughout this section, $G$  denotes a range finite directed
graph with no sources; $\sG$ denotes the groupoid associated with $G$;
and $\sA$ is the \cstar algebra constructed from $G$ or $\sG$. 
We shall need to use the convenient fact that elements of $\sA$ can be
identified with continuous functions on $\sG$ which vanish at
infinity.  This is a consequence of the range finitness of $G$, 
which implies that the groupoid $\sG$ is amenable
(\cite[Corollary~5.5]{MR98g:46083} and~\cite[Theorem~4.2]{MR1962477}). 
 Amenability, in turn, implies
that $C^*(\sG) = C^*_{\text{red}}(\sG)$ \cite[p.~92]{MR82h:46075}.  
Finally, Proposition II.4.2 in~\cite{MR82h:46075} allows us to
identify the elements of $C^*(\sG)$ with (some of the) elements in
$C_0(\sG)$, the continuous functions on $\sG$ vanishing at infinity.

\begin{definition} \label{d:spectrum}
Let $\sB \subseteq \sA$ be a bimodule over $\sD$.  Define the
\emph{spectrum} of $\sB$ to be 
\begin{displaymath}
  \sigma(\sB) = \{(x,k,y) \in \sG \mid \text{there is } f \in \sB \text{
    such that } f(x,k,y) \neq 0 \}.
\end{displaymath}
For any open subset  $P$  of $\sG$, we let
\begin{displaymath}
  A(P) = \{ f \in \sA \mid f(x,k,y) = 0 \text{ for all } (x,k,y )
  \notin P \}.
\end{displaymath}
\end{definition}

It is easy to check that $A(P)$, which consists of all 
functions in $\sA$ which
are supported on $P$, is a bimodule over $\sD$ and that
$\sigma(A(P)) = P$.  It is also trivial to see that 
 $\sB \subseteq A(\sigma(\sB))$.  But it is not 
always true that $\sB = A(\sigma(\sB))$; a
counterexample is given in~\cite{hop_peters}.  
 Thus, the spectral theorem for bimodules for
graph \cstar algebras differs from the theorem for
groupoid \cstar algebras where the groupoid is $r$-discrete, amenable
and principal~\cite{MR90m:46098}.  

In fact, the existence of a counterexample depends exactly on the
presence of a loop in the graph.

\begin{proposition} \label{p:counterex}
Let $G$ be a range finite directed graph with no sources.  There is in
$\sA$ a bimodule $\sB$ over $\sD$ such that
$\sB \ne A(\sigma(\sB))$ if, and only if, $G$ has a loop.
\end{proposition}

\begin{proof}
Suppose $G$ has a loop, say
$\alpha = \alpha_1 \alpha_2 \dots \alpha_k$ with
$s(\alpha_k) = s(\alpha) = r(\alpha) = r(\alpha_1)$.  Let
$x = \alpha \alpha \alpha \dots $ in $P(G)$.  Write
$S_{\alpha} = S_{\alpha_1} \dots S_{\alpha_k}$,
$v = r(\alpha)$, and $\Phi = P_v + S_{\alpha}$.  Let
$\sB$ denote the bimodule generated by $\Phi$; this bimodule is the
norm closure of all finite sums
$\sum f_i \Phi g_i$, with $f_i, g_i \in \sD$.

Now $\sigma(\Phi) = \sigma(P_v) \cup \sigma(S_{\alpha})$ and
\begin{align*}
\sigma(P_v) &= \{(z,0,z) \mid z \in P(G), r(z) = v\}, \\
\sigma(S_{\alpha}) &= \{ (\alpha z, k, z) \mid
z \in P(G), r(z) = s(\alpha) = s(\alpha_k) \}.
\end{align*}
Since $(x,0,x) \in \sigma(P_v)$ and
$(x,k,x) \in \sigma(S_{\alpha})$, $(x,0,x)$ and $(x,k,x)$ both lie in
$\sigma(\Phi)$. 

Viewing elements of $\sD$ as functions on the groupoid supported on the 
unit space and using 
$\Phi = \chi_{\sigma(P_v)} + \chi_{\sigma(S_{\alpha})} =
\chi_{\sigma(P_v) \cup \sigma(S_{\alpha})}$, we have
\begin{align*}
f\Phi (x,0,x) &= f(x,0,x) \Phi(x,0,x) = f(x,0,x), \\
f\Phi (x,k,x) &= f(x,0,x) \Phi(x,k,x) = f(x,0,x).
\end{align*}
Similarly, $\Phi f(x,0,x) = f(x,0,x) = \Phi f (x,k,x)$.

So, for any $f,g \in \sD$, we have
$f \Phi g (x,0,x) = f \Phi g (x,k,x)$.  The same equality is valid for 
sums and extends to the norm closure: if $h \in \sB$, then
$h(x,0,x) = h(x,k,x)$.  Since there are elements
$h' \in A(\sigma(\sB))$ with $h'(x,0,x) \ne h'(x,k,x)$, this shows
that $\sB  \ne A(\sigma(\sB))$.

On the other hand, if $G$ has no loops, $C^*(G)$ is AF and
$\sG$ is a principal grooupoid, so all bimodules satisfy
$\sB = A(\sigma(\sB))$ by the Muhly-Solel spectral theorem for
bimodules~\cite{MR90m:46098}. 
\end{proof}

The spectral theorem for bimodules below
provides two necessary and 
sufficient conditions for a bimodule to be
determined by it spectrum.  These are, in fact, the two equivalent
conditions in Theorem~\ref{t:stb_prelim}.  A third equivalent
condition will be given in section~\ref{s:stb3}.

\begin{theorem}[Spectral Theorem for Bimodules] \label{STB}
Let $G$ be a range finite directed graph with no sources. 
 Let $\sB \subseteq \sA$ be a bimodule over $\sD$.
Then the following statements are equivalent\/{\rm :}
\begin{enumerate}
\item $\sB = A(\sigma(\sB))$. \label{c:detspec}
\item $\sB$ is generated by the Cuntz-Krieger partial isometries which
  it contains.  \label{c:CKgen}
\item $\sB$ is invariant under the gauge automorphisms. \label{c:gaugeinv}
\end{enumerate}
\end{theorem}

\begin{proof}
Since the equivalence of~\ref{c:CKgen} and~\ref{c:gaugeinv} has
already been established in Theorem~\ref{t:stb_prelim}, it suffices to prove
the equivalence of~\ref{c:detspec} and~\ref{c:CKgen}.

To show that $\ref{c:detspec} \Rightarrow \ref{c:CKgen}$ we need
to show that whenever $P$ is an open subset of $\sG$, $A(P)$ is
generated by the Cuntz-Krieger partial isometries which it contains.
Let $\sB$ be the bimodule which is generated by the Cuntz-Krieger
partial isometries which are in $A(P)$.  Clearly,
$\sB \subseteq A(P)$; we must show the reverse containment.

Assume that $\SSadj{\alpha}{\beta} \in A(P)$.  So
$Z(\alpha,\beta) \subseteq P$.  Let $f$ be a continuous function with
support in $Z(\alpha, \beta)$.  Define 
$g \colon Z(\alpha,\alpha) \rightarrow \bbC$ by
\begin{displaymath}
 g(\alpha\gamma,0,\alpha\gamma) = f(\alpha\gamma, |\alpha|-|\beta|,
\beta\gamma). 
\end{displaymath}
Then,
\begin{align*}
g\cdot \SSadj{\alpha}{\beta} (x,k,y) &= g(x,0,x)
\SSadj{\alpha}{\beta}(x,k,y) \\
&=
\begin{cases}
g(x,o,x), \qquad &\text{if } x = \alpha \gamma, k = |\beta| - |\alpha|,
y = \beta \gamma \\
0,  &\text{otherwise}
\end{cases} \\
&=
\begin{cases}
f(x,k,y), \qquad &\text{if } x = \alpha \gamma, k = |\beta| - |\alpha|,
y = \beta \gamma \\
0, &\text{otherwise}
\end{cases} \\
&= f(x,k,y).
\end{align*}
Since $g \in \sD$ and $\SSadj{\alpha}{\beta} \in \sB$, we have
$f \in \sB$.  

Thus, $\sB$ contains any continuous function supported on a set of the
form $Z(\alpha,\beta) \subseteq P$.  Any compact open subset of $P$
can be written as a finite union of sets of the form 
$Z(\alpha,\beta)$, so $\sB$ contains any $f \in \sA$ which is
supported on a compact open subset of $P$.  But any compact subset of
$P$ is contained in a compact open subset of $P$, so $\sB$ contains all
$f$ which are supported on a compact subset of $P$.  These functions
are dense in $A(P)$ in the \cstar norm, so $A(P) \subseteq \sB$.

To prove that $\ref{c:CKgen} \Rightarrow \ref{c:detspec}$, assume that
$\sB$ is a bimodule over $\sD$ which is generated by the Cuntz-Krieger
partial isometries which it contains.
We first show that $\ds \sigma(\sB) = \cup Z(\alpha,\beta)$, where the
union is taken over all $\alpha,\beta$ such that
$\SSadj{\alpha}{\beta} \in \sB$.  Indeed, if $p$ is in this union, then
there is $\SSadj{\alpha}{\beta} \in \sB$ such that
$\SSadj{\alpha}{\beta} (p) = 1$, so $p \in \sigma(\sB)$.  On the other
hand, if $p$ is not in the union, then $\SSadj{\alpha}{\beta} (p) = 0$
for all $\SSadj{\alpha}{\beta} \in \sB$.  Write $p = (x,k,y)$ and let
$f$ and $g$ be in $\sD$.  Then
\begin{displaymath}
 f \cdot  \SSadj{\alpha}{\beta} \cdot g (p) = 
f(x,0,x) \SSadj{\alpha}{\beta}(x,k,y) g(y,0,y) = 0.
\end{displaymath}
It follows that all elements of the bimodule generated by the Cuntz-Krieger
partial isometries in $\sB$ vanish at $p$.  But the Cuntz-Krieger
partial isometries generate $\sB$ itself, so $p \notin \sigma(\sB)$.  

We already know that $\sB \subseteq A(\sigma(\sB))$; to show the reverse
containment it is sufficient (since any $A(P)$ is generated by the
Cuntz-Krieger partial isometries which it contains)
to show that if 
$\SSadj{\alpha}{\beta} \in A(\sigma(\sB))$ then 
$\SSadj{\alpha}{\beta} \in \sB$.
 
Let $\SSadj{\alpha}{\beta} \in A(\sigma(\sB))$, so that 
$Z(\alpha,\beta) \subseteq \sigma(\sB)$.  Let $p \in Z(\alpha,\beta)$.
Then  there is $\SSadj {\nu} {\mu} \in \sB$ such
that $p \in Z(\nu, \mu) \subseteq \sigma(\sB)$.  Since
$Z(\alpha,\beta) \cap Z(\nu,\mu) \ne \emptyset$, we have either
$Z(\alpha,\beta) \subseteq Z(\nu, \mu)$ or
$Z((\nu, \mu) \subseteq Z(\alpha,\beta)$.  

If $Z(\alpha,\beta) \subseteq Z(\nu, \mu)$, then there is a finite
path $\epsilon$ such that $\alpha = \nu \epsilon$ and
$\beta = \mu \epsilon$.  A routine calculation shows that
$\SSadj {\alpha}{\beta} = \SSadj {\alpha}{\alpha} \SSadj {\nu} {\mu}
\SSadj {\beta} {\beta}$.  But $\SSadj {\alpha}{\alpha}$ and
$\SSadj {\beta} {\beta}$ are in $\sD$ and 
$\SSadj {\nu} {\mu} \in \sB$, so $\SSadj{\alpha}{\beta} \in \sB$, as
desired.

Suppose, on the other hand, that for any point $p \in Z(\alpha, \beta)$,
the set $Z(\nu, \mu)$ obtained as above is a subset of 
$Z(\alpha,\beta)$.  Then these sets form an open 
cover for $Z(\alpha, \beta)$.  Since $Z(\alpha,\beta)$ is compact, we
can find a finite subcover.  It is routine to arrange that this
subcover is disjoint (without losing the property that the associated
Cuntz-Krieger partial isometries are in $\sB$).  Thus, we can write
$\ds Z(\alpha, \beta) = \cup^{n}_{i=1} Z(\nu_i, \mu_i)$, 
a finite disjoint
union with all $\SSadj {\nu_i}{\mu_i} \in \sB$.  But then
$\ds \SSadj {\alpha}{\beta} = \sum^n_{i=1} \SSadj {\nu_i}{\mu_i}$ and so
$\SSadj {\alpha} {\beta} \in \sB$.
\end{proof}

\section{Cocycles} \label{s:cocycles}

As usual $G$ is a range finite directed graph with no sources.
Most of our attention will be focused on cocycles
defined on the associated Cuntz-Krieger
groupoid $\sG$.
 
A real valued $1$-cocycle on $\sG$ is a
continuous function $c \colon \sG \to \bbR$ which satisfies
the cocycle condition
$c(x, k, y) + c(y, l, z) = c(x, k+l, z)$, for all composable pairs. 
It follows that $c(x, 0, x) = 0$, for $x$ in path space $P$, and that
$c((x, k, y)^{-1}) = c(y, -k, x) = -c(x, k, y)$, 
for all $(x, k, y) \in \sG$. 
The set of all $1$-cocylces forms a group 
under addition, denoted by $Z^1(\sG, \bbR)$.

A simple example of a cocycle is the one given by the formula
$c(x,k,y)=k$.  This cocycle is intimately related to the gauge
automorphisms: for any $f \in C_c(\sG)$,
$\gamma_z(f)(x,k,y) = z^k f(x,k,y)$.  More generally,
any cocycle $c$ gives rise to the one-paramenter automorphism group
\[ 
\eta_z(f)(x, k, y) = z^{c(x, k, y)} f(x, k, y) .
\]
Each $\eta_z$ is a $*$-automorphism of $C_c(\sG)$ onto itself; it is 
not hard to show that this automorphism preserves
the \cstar norm and so extends to an automorphism of $\sA$ with
the formula above.  

\begin{remark} \label{r:Hinfty}
For each point
$(x, k, y) \in \sG$, the map $f \mapsto f(x, k, y)$ 
on $C_c(\sG)$ is
 decreasing with respect to the $||\cdot ||_{\infty}$ norm, and hence
also decreasing with respect to the \cstar norm
(\cite[Prop.~II.4.1]{MR82h:46075}).  
So these maps extend to  continuous linear functionals on 
$\sA$. 
 If $f \in \sA$, we consider all 
functions of the form $t \to \rho(\eta_t(f))$, where $\rho$
is a linear functional of the type above. 
Given $f \in \sA$, it is easy to check that
$t \to \rho(\eta_t(f))$ is an $H^{\infty}$-function on $\bbR$ for all 
linear functionals of this form if, and only if,
$f$ is supported on $\{ (x, k, y) \in \sG \mid c(x, k, y) \geq 0\}$.
(Note: when we write $\eta_t$ it refers, of course, to
$\eta_z$ for $z=e^{it}$; use of the real variable is more appropriate
when discussing $H^{\infty}$.)

To emphasize the connection with analyticity, consider the simplest
possible graph: the graph $G$ consisting of a single vertex and a
single loop.  The associated groupoid for $G$ is the group
of integers, $\bbZ$.  The groupoid \cstar algebra (well, really, the
group \cstar algebra) is isomorphic to $C(\bbT)$.  Briefly,
$C_c(\bbZ)$ is a \star algebra in which   the
multiplication is convolution
and the \cstar norm of a function $f \in C_c(\bbZ)$ is the 
$\|\ \|_{\infty}$ norm of the funtion
$\theta \mapsto \sum_{n} f(n) e^{in\theta}$, 
$\theta \in \bbT$. Therefore, $C^*(\bbZ)$ is identified with the $C_0$
functions on $\bbZ$ which are the Fourier coefficients of functions
in $C(\bbT)$.

A cocycle on $\bbZ$ is determined by its value at $1$; so the only
one of interest is $c(n) = n$.  Let $\alpha_t$ be the associated one
parameter family of automorphisms acting on $C^*(\bbZ)$.  Since
$\alpha_t(f)(n) = e^{int}f(n)$, when transferred via the inverse
Fourier transform to $C(\bbT)$, the automorphism group acts by
translation: $\alpha_t(\phi)(\theta) = \phi(\theta + t)$.

For each $n \in \bbZ$, let $\rho_n$ be the linear functional on
$C^*(\bbZ)$ given by $f \mapsto f(n)$.  Transferred to $C(\bbT)$, this
is $\phi \mapsto \int \phi(\theta) e^{in \theta}\,d\theta$, where
$d\theta$ is normalized Lebesgue measure on $\bbT$.  The closed linear
span of the functionals of this type can be identified with the
complex valued measures on $\bbT$ which are absolutely continuous with
respect to Lebesgue measure; i.e. with $L^1(\bbT, d\theta)$.  Thus,
``evaluation functionals'' do not span the dual space of
$C(\bbT)$.

All the same, if $\phi \in C(\bbT)$ and 
$t \mapsto \rho_n(\alpha_t(\phi))$  is in $H^{\infty}(\bbR)$
for all $n$, then
$t \mapsto \rho(\alpha_t(\phi))$ is in $H^{\infty}(\bbR)$ for all 
$\rho \in C(\bbT)^*$.  This happens exactly when the Fourier
transform $\hat{\phi}$ is supported on $\bbZ^+ = c^{-1}[0,\infty)$.
Thus the analytic algebra associated with the cocycle $c$ is just the
disk algebra, $A(\bbD)$.
\end{remark}

\begin{definition} \label{d:shift}
Let $\sS \colon P \to P$ be the shift map; thus 
if $x = x_1 x_2 \dots$ is a path with terminal
edge $x_1$, $\sS(x) = x_2 x_3 \dots$.
\end{definition}

Note that $\sS$ is a continuous map, in fact it is a local homeomorphism.

Let $\sG(k)$ denote 
$\{ (x, l, y) \in \sG \mid l = k\}$ and
 $C(P)$  denote the space of continuous functions
on $P$.
 We now give an
example of a class of cocycles.

\begin{example} \label{e:cocy}
 Let $f \in C(P)$ and define $c$ on $\sG(k) \ (k > 0)$ by
\[ 
c(x, k, y) = \sum_{j=0}^{k-1} f(\sS^jx) + \sum_{j=k}^{\infty}
[f(\sS^jx) - f(\sS^{j-k}y)] 
\]
for $(x, k, y) \in \sG(k)$.  Observe that the 
infinite sum has only finitely many nonzero terms.
For $k = 0$, set 
$ c(x, 0, y) = \sum_{j=0}^{\infty} [f(\sS^jx) - f(\sS^jy)]$, 
and for $k$ negative  set
$c(x, k, y) = -c(y, -k, x)$.

To verify the cocycle condition, let 
$(x, k, y)$, $(y, l, z) \in \sG$ with $k, l \geq 0$. Then
\begin{align*}
c(x, k, y) + c(y, l, z) &= \sum_{j=0}^{k-1} f(\sS^jx) + 
\sum_{j=k}^{\infty} [f(\sS^jx) - f(\sS^{j-k}y)] \\
                        &+ \sum_{i=0}^{l-1} f(\sS^iy) + 
\sum_{i=l}^{\infty} [f(\sS^iy) - f(\sS^{i-l}z)] \\
                                &= \sum_{j=0}^{k+l-1} f(\sS^jx) + 
\sum_{j=k+l}^{\infty} [f(\sS^jx) - f(\sS^{j-k-l}z)] \\
                                &= c(x, k+l, z).
\end{align*}
The other cases are similar. Finally, observe that if $f$ is 
continuous on $P$, then the cocycle $c$ is
continuous in the topology on $\sG$.
 \end{example}

Note that the cocycle $c(x, k, z) = k$, which generates the gauge automorphisms,
is produced by the constant function $f(x) = 1$. 

\begin{theorem} \label{t:cocy} 
Let $G$ be a range finite directed graph with no sources.
Let $\sG$ be the associated Cuntz-Krieger groupoid and 
$P$ the path space of $G$,
identified with the unit space $\sG^0$.
Then there is a bijection 
$C(P) \longleftrightarrow Z^1(\sG,\bbR)$ given
 as follows\/{\rm :} for $f \in C(P)$, let
$c_f$ denote the cocycle constructed in Example~\ref{e:cocy}. 
For $c \in Z^1$, let $f_c(x) = c(x, 1, \sS x)$.
Then the two maps are inverses of each 
other\/{\rm :} $f = f_{c_f}$ and $c = c_{f_c}$.
\end{theorem}

\begin{proof} Let $c \in Z^1(\sG, \bbR) $ be given, and $(x, 0, y) \in
  \sG(0)$.
 Note that $(\sS{x}, 0, \sS{y})$
also belongs  to $\sG(0)$. From the cocycle condition we have
\begin{equation} \label{e:cocy1}
  \ c(x, 1, \sS{y}) = c(x, 1, \sS{x}) + c(\sS{x}, 0, \sS{y}) =
 c(x, 0, y) + c(y, 1, \sS{y}) .
\end{equation}
With $f(x) = f_c(x) = c(x, 1, \sS{x})$, $f \in C(P)$, and we can
rewrite equation~\eqref{e:cocy1} as 
\begin{equation} \label{e:cocy2}
f(x) - f(y) = c(x, 0, y) - c(\sS{x}, 0, \sS{y}).
\end{equation}
Replacing $x, y$ with $\sS^j{x}, \sS^j{y}$ and summing over $j= 0, 1, \dots$ we obtain
\[ c(x, 0, y) = \sum_{j=0}^{\infty} [f(\sS^j{x}) - f(\sS^j{y})] .\]
Since for sufficiently large $j,\ \sS^j{x} = \sS^j{y}$, the sum above is actually finite.

Now let $(x, k, y) \in \sG(k)$ with $k > 0$.  Then $(\sS^k{x}, 0, y) \in \sG(0)$, and we have
\begin{align*}
 c(x, k, y) &= \sum_{j=0}^{k-1} c(\sS^{j}{x}, 1, \sS^{j+1}{x}) + c(\sS^k{x}, 0, y)\\
                &= \sum_{j=0}^{k-1} f(\sS^j{x}) + \sum_{j=k}^{\infty} 
             [f(\sS^j{x}) - f(\sS^{j-k}{y})] \\
                &= c_f(x, k, y).
\end{align*}
The case $k < 0$ is similar.

Conversely, given $f \in C(P)$, define the cocycle $c$ as in Example~\ref{e:cocy}. But then
$c(x, 1, \sS{x}) = f(x)$, so $f = f_{c_f}$.
\end{proof}

\begin{remark} 
If $G$ is a finite directed graph, 
then $\sG(0)$ is the AF-groupoid associated with the
stationary Bratteli diagram which at level $n$ has a copy of 
the vertices $V$, and admits an edge from vertex $v$
at level $n$ to vertex $w$ at level $n+1$ if, and
only if, $G$ has an edge 
from $v$ to $w$. The restriction of a cocycle $c$
on the groupoid $\sG$ to the subgroupoid $\sG(0)$ gives a cocycle
 on the AF groupoid $\sG(0)$.  This class of cocycles
has not been systematically studied. It is, however, a proper 
subclass of $Z^1(\sG(0), \bbR)$, as we shall see in
the subsection on integer-valued cocycles.
\end{remark}

\begin{proposition} \label{p:unif_conv}
Let $f$ and $f_n \ (n \in \bbN)$ be continuous functions 
on $P$.  Then $f_n \to f$ uniformly on compact
subsets of $P$ if, and only if, $c_{f_n} \to c_f$ uniformly on 
compact subsets of $\sG$.
\end{proposition}

\begin{proof} Suppose $f_n \to f$ uniformly on each of the sets 
$Z(\alpha)$, $\alpha \in F$. Given $\alpha, \beta \in F$ 
and the basic  compact open set
 $Z(\alpha, \beta) \subset \sG(k)$ where $k = |\alpha| - |\beta| \geq 0$,
we have, for $(x, k, y) \in Z(\alpha, \beta)$,
\begin{multline*}
c_{f_n}(x, k, y) - c_{f}(x, k, y) = \\
\sum_{j=0}^{k-1} [f_n(\sS^j{x}) - f(\sS^j{x})] + 
\sum_{j=k}^{|\alpha|} ([f_n(\sS^j{x}) - f(\sS^j{x})]
- [f_n(\sS^{j-k}{y}) - f(\sS^{j-k}{y})])
\end{multline*}
and this converges uniformly to zero on $Z(\alpha, \beta)$.  A similar
argument applies
 when $|\alpha|- |\beta| < 0$.

For the converse, suppose that $Z(\alpha) \subset P$.
Write $\alpha = x_1 x_2 \dots x_n$, and let $\beta$ be the empty
string if $|\alpha|=1$ and
 $\beta = x_2 x_3 \dots x_n$ otherwise.
Note that $x \in Z(\alpha)$ if, and only if,
$(x, 1, \sS{x}) \in Z(\alpha, \beta)$. Since $f(x) = c_f(x, 1, \sS{x})$ 
and  $c_{f_n}$ converges uniformly
to $c_f$ on $Z(\alpha, \beta)$, it follows that $f_n$ converges 
uniformly to $f$ on $Z(\alpha)$.
\end{proof}

\begin{proposition} \label{p:loc_const}
Let $f \in C(P)$ and let $c_f$ be the corresponding cocycle
  on $\sG$. 
 Then $f$ is locally constant on $P$
if, and only if, $c_f$ is locally constant on $\sG$.
\end{proposition}

\begin{proof} Assume $c_f$ is locally constant. Given $x \in P$, 
let $Z(\alpha, \beta)$ be a neighborhood of
$(x, 1, \sS{x})$ on which $c_f$ is constant. Since 
$f(u) = c_f(u, 1, \sS{u})$ for all  $u \in P$, it follows that $f$ is constant
on $Z(\alpha)$.

Suppose now that $f$ is locally constant, and let $(x, k, y) \in \sG$
be given. 
 We suppose that $k \geq 0$ (the case
$k < 0$ is analogous).  There is
$n \geq k$ such that $\sS^j{x} = \sS^{j-k}y$ for all $j \geq n$.  
Therefore we can write
\[ 
c_f(x, k, y) = \sum_{j=0}^{k-1} f(\sS^j{x}) + \sum_{j=k}^n
[f(\sS^j{x}) - f(\sS^{j-k}{y})] .
\]
For $p$ chosen sufficiently large, if we let
$\alpha = x_1\dots x_p$ and $\beta = y_1 \dots y_{p-k}$, then
$x_i = y_{i-k}$, for all $i \geq p+1$;
$Z(\sS^j(\alpha))$ is a clopen neighborhood of 
$\sS^j(x)$ on which $f$ is constant, for $j=1,\dots,n$; and
$Z(\sS^j(\beta))$ is a clopen neighborhood of
$\sS^j(y)$ on which $f$ is constant, for $j=1,\dots,n-k$.
Then $(x, k, y) \in Z(\alpha, \beta)$ and
$c_f$ is constant on $Z(\alpha, \beta)$.
\end{proof}

\begin{remark}
Proposition~\ref{p:loc_const} applies, in particular, whenever $f$ has
 finite range.
\end{remark}

\begin{definition} \label{d:Z_0}
Let $Z^1_0(\sG, \bbR)$ denote the subset of
  $Z^1(\sG, \bbR)$ consisting of those cocycles $c$ which
vanish precisely on the unit space  $\sG^0$.
\end{definition}

\begin{remark} Every cocycle in $Z^1(\sG, \bbR)$ necessarily vanishes
on  $\sG^0$. Also, note that
$Z^1_0$ is not a subgroup of $Z^1$; indeed, $ 0 \notin Z^1_0$.
\end{remark}

\subsection{Bounded Cocycles}
In the context of AF algebras and their groupoids, bounded 
cocycles are of special interest due to the connection
between bounded cocycles and reflexive subalgebras of
 AF algebras (cf.~\cite{MR95c:46091}).
Thus, it is natural to
investigate the role of bounded cocycles on Cuntz-Krieger groupoids.

A point $x \in P$
is \emph{periodic} if $\sS^k{x} = x$ for some $k > 0$. Note that 
the existence of a periodic point in $P$ is
equivalent to the existence of a loop in the graph $G$. Recall from 
\cite{MR99i:46049} that $G$ has no loops if, and only if,
C$^*(G)$ is an AF-algebra.

\begin{proposition}  \label{p:nobdd_cocy}
Let $G$ be a range finite directed graph with no sources.
 Then $Z^1_0(\sG, \bbR)$ contains no bounded cocycle if, and only if,
$G$ contains a loop.
\footnote{The implication $Z_0^1(\sG,\bbR)$
  contains no bounded cocycle $\Rightarrow$ $G$ contains a loop is due 
  to Allan Donsig.  The authors thank Donsig for giving permission to
  use this result here.} 
\end{proposition}

\begin{proof} 
First, assume that $G$ contains a loop.
Then there is a periodic point, say
 $x = x_1 \dots x_p x_1 \dots x_p x_1 \dots $, in $P$.
Let $c \in Z^1_0(\sG, \bbR)$ and
  let $f(y) = c(y, 1, \sS{y}), \ y \in P$.
For $k \geq 1$, we have
\begin{align*}
c(x, kp, \sS^{kp}{x}) &= f(x) + f(\sS{x}) + \cdots + f(\sS^{p-1}x) \\
                          &+ f(\sS^p{x}) + f(\sS(\sS^p{x})) + 
                           \cdots + f(\sS^{p-1}(\sS^p{x})) \\
                          &+ \cdots \\
                          &+ f(\sS^{(k-1)p}{x}) +
                          f(\sS^{(k-1)p}(\sS{x})) 
                           + \cdots + f(\sS^{(k-1)p}(\sS^{p-1}{x})).
\end{align*}
Using $\sS^p{x} = x$, this reduces to
\[
 c(x,kp,\sS^{kp}{x}) = k[f(x) + f(\sS{x}) + \cdots + f(\sS^{p-1}{x})], 
\quad \text{for all } k \geq 1.
\]
But $\sum_{j=0}^{p-1} f(\sS^j{x}) = c(x, p, \sS^p{x}) \neq 0 $ since 
$(x, p, \sS^p{x}) = (x, p, x) \notin \sG_0$.
As $k$ is arbitrary, $c$ is unbounded.

Now assume that $G$ contains no loop.  Let
$a \colon G \to \bbR^+$ be a function with the property that for each
edge $e$,
$T_e = \{ f \in G \mid a(f) > a(e) \}$ is finite and
\[
a(e) > \sum_{f \notin T_e} a(f).
\]
(This is easily done after $G$ is arranged as a sequence.)

Now define a continuous, locally constant function
$f \colon C(P) \to \bbR$ by $f(x) = a(x_1)$, for $x=x_1x_2\dots$.  
In other words, $f$ has
the value $a(x_1)$ on $Z(x_1) \subseteq C(P)$.  Let $c$ be the cocycle 
associated with $f$, as in example~\ref{e:cocy}. Clearly, $c$ is
bounded.  Since any cocycle
vanishes on the unit space, we just need to show that $c$ is non-zero
off the unit space.

Because $G$ has no loops, there are no points of the form
$(x,k,x) \in \sG$ with $k \ne 0$.  Therefore, we just need to show
that $c(x,k,y) \ne 0$ whenever  $x \ne y$.
  Thanks to the cocycle
property, we may without loss of generality assume that $k \geq 0$.  
Write
\begin{align*}
x &= x_1 \dots x_{p+k} z_1 z_2 \dots , \\
y &= y_1 \dots y_p z_1 z_2 \dots .
\end{align*}
We then have 
\begin{align*}
c(x,k,y) &= \sum_{j=0}^{k-1} f(\sS^j x) + \sum_{j=k} ^{\infty}
 [f(\sS^j x) -f(\sS^{j-k} y)] \\
&= \sum_{j=1}^{p+k} a(x_j) - \sum_{j=1}^p a(y_j).
\end{align*}
Since $G$ has no loops, a given edge $e$ may appear at most once in
each of the paths $x_1 \dots x_{p+k}$ and $y_1 \dots y_p$.  Some edges 
may appear in both paths, but then the terms cancel.  
Since $x \ne y$, there is an edge $e$ 
amongst $x_1, \dots, x_{p+k}, y_1, \dots, y_p$ which
appears once only and for which $a(e)$ is maximal;  the summation
property for $a$ then guarantees that $c(x,k,y) \ne 0$.
\end{proof}

Recall that a graph is \emph{transitive} if there is a path from any
vertex to any other vertex.
If a directed graph $G$ is finite and transitive, it 
satisfies Cuntz and Krieger's condition for $C^*(G)$ to be simple.

\begin{proposition} \label{p:bdd_cocy}
Let $G$ be a finite, transitive, directed graph, and let $c$ 
be a bounded cocycle on the AF groupoid $\sG(0)$. Suppose
that C$^*(\sG(0))$ is simple. Then $c$ extends to a cocyle on $\sG$.  
Furthermore, if $c$ vanishes precisely on
the unit space $\sG^0$, then the extension can be chosen to
 vanish precisely on $\sG^0$.
\end{proposition}

\begin{proof}  By~\cite[p.~112]{MR82h:46075},
since C$^*(\sG(0))$ is simple and $c$ is bounded, $c$ is a coboundary: that is,
there is a continuous function $b$ on $P$ 
so that $c(x, 0, y) = b(x) - b(y)$.

Choose a point $x_0 \in P$, and let $[x_0]$ denote the equivalence
class 
of $x_0$ in $\sG(0)$: $[x_0] = \{ y \in P \mid (y, 0, x_0) \in \sG(0) \}$. 
Since C$^*(\sG(0))$ is
simple, it follows from~\cite{MR82h:46075}
that the equivalence class of any point is dense; thus $[x_0]$ is dense.

We shall construct $f\in C(P)$ such
that $c_f$ extends $c$.
 Begin by setting $f(x_0) = 0$.
For $x \in [x_0]$, define
\[
 f(x) = c(x, 0, x_0) - c(\sS{x}, 0, \sS{x_0}) .
\]
  The cocycle property for $c$ shows that if $x, y \in [x_0]$,  then
\begin{equation} \label{e:cocy3}
 f(x) - f(y) = c(x, 0, y) - c(\sS{x}, 0, \sS{y}). 
\end{equation}
Since $P$ is compact (because $G$ is finite), Hausdorff and first
countable, it is metrizable.  Let $\rho$ be a metric for $P$. 
Since a continuous function on a compact metric
space is uniformly continuous, both $b$ and $b\circ \sS$ are uniformly
continuous.  Thus, given $\epsilon > 0 $
there is a clopen cover $\{ Z(\alpha) \mid \alpha \in \bbA \}$ with 
$\bbA$ finite and such that for 
$x, y \in Z(\alpha)$, 
both $|b(x) - b(y)| < \frac{\epsilon}{2}$ and 
$|b\circ \sS(x) - b\circ \sS(y)| < \frac{\epsilon}{2}$.

Now for $x, y \in [x_0] \cap Z(\alpha)$ we have
\begin{align*}
|f(x) - f(y)| &= |c(x, 0, y) - c(\sS{x}, 0, \sS{y})| \\
                  &= |(b(x) - b(y)) - (b\circ \sS(x) - b\circ \sS(y))| \\
                  &\leq |b(x) - b(y)| + |b\circ \sS(x) - b\circ \sS(y)| \\
                  &< \frac{\epsilon}{2} + \frac{\epsilon}{2} = \epsilon.
\end{align*}

As $f$ is uniformly continuous on a dense subset, it admits a
continuous extension to $P$, also denoted by $f$. Note
that since $\{ (x, 0, y) \mid x, y \in [x_0] \}$ is dense in $\sG(0)$, it
follows that equation~\eqref{e:cocy3} holds for
$x, y \in P$. By the same argument used in the proof of Theorem~\ref{t:cocy},
\[ c(x, 0, y) = \sum_{j=0}^{\infty} [f(\sS^j{x}) - f(\sS^j{y})] .\]
It is now immediate that the cocycle $c_f$ (see Example~\ref{e:cocy}) extends $c$. 

Fix $k \in \bbZ$; we claim that $c_f$ is bounded on each $\sG(k)$. 
Of course if $k = 0$ this holds by assumption. If
$k > 0$ and $(x, k, y) \in \sG(k)$ we can write
\begin{align*}
 c_f(x, k, y) &=\sum_{j=0}^{k-1} f(\sS^j{x}) + \sum_{j=k}^{\infty} 
           [f(\sS^j{x}) - f(\sS^{j-k}{y})] \\
                 & =\sum_{j=0}^{k-1} f(\sS^j{x}) + c(\sS^k{x}, 0, y)
\end{align*}
Since $f$ is bounded, $k$ is fixed, and $c$ is bounded on $\sG(0)$, 
$c$ is bounded on $\sG(k)$. From the
cocycle property it follows that $c$ is bounded on $\sG(k)$, 
 for $k$ negative as well.

Finally, suppose that $c$ vanishes precisely on $\sG^0$. If $L$ is a
constant, $g \in C(P)$, and $g = f + L$, the cocycle $c_g$ is related to $c_f$ 
by $c_g (x, k, y) = c_f(x, k, y) + kL$, for $(x, k, y) \in \sG(k)$.
In particular, $c_f$ and $c_g$ agree on $\sG(0)$, and 
$c_g(x, 1, y) = c_f(x, 1, y) + L$, for $(x, 1, y) \in \sG(1)$.
Since $c_f$ is bounded on $\sG(1)$, we can pick $L$ sufficiently 
large so that $c_g$ is positive on $\sG(1)$.
For $k > 1$, any element of $\sG(k)$ is a product of $k$ elements in 
$\sG(1)$, so that $c_g$ is positive on $\sG(k)$.
For $k < 0$, the cocycle property guarantees that $c_g$ is negatve on $\sG(k)$. 
Hence, the equation $c_g = 0$
has exactly the same solutions as $c = 0$ on $\sG(0)$, namely $c_g =0$ 
precisely on the unit space $\sG^0$.
\end{proof}

\begin{remark}\label{r:cuntz}
In case $G$ is the directed graph with a single vertex and $n$ 
edges (i.e., loop edges), $\sG$ is the
Cuntz groupoid and $\sG(0)$ is the UHF$(n^{\infty})$ groupoid.
 In this case the hypotheses of  Proposition~\ref{p:bdd_cocy}
are satisfied,
and hence any bounded cocycle which vanishes precisely on the unit 
space extends to a cocycle on $\sG$ vanishing
precisely on $\sG^0$.  In particular, this applies to the refinment
cocycle on the UHF$(n^{\infty})$ groupoid.
\end{remark}
 
\subsection{Integer-valued Cocycles}
Among the most studied classes of cocycles on AF groupoids are the 
integer-valued cocycles, due to their connection
with dynamical systems.  Thus it is natural to examine integer-valued 
cocycles on Cuntz-Krieger groupoids.

\begin{lemma} \label{l:integer} 
Let $\sG$ be a Cuntz-Krieger groupoid, and suppose $\sG(0)$ has a dense equivalence
class. \/{\rm (}In particular, that will be the case 
when C$^*(\sG(0))$ is simple.\/{\rm )} 
Suppose $c$ is an integer valued cocycle
defined on $\sG(0)$. Then 
\begin{enumerate}
\item
$c$ admits an extension to a cocycle on  $\sG$ if,
 and only if, $c$ admits an 
extension to an integer valued cocycle on $\sG$.
\item
If  $c$ vanishes precisely on the unit space $\sG^0$,  
then  $c$ admits an extension to 
$\sG$  with this property if, and only if,  $c$ admits an 
integer valued extension to $\sG$ vanishing
precisely on $\sG^0$.
\end{enumerate}
\end{lemma}

\begin{proof}  By Theorem~\ref{t:cocy} any extension of $c$ to $\sG$ 
is of the form $c_f$, for some $f \in C(P)$.
Furthermore, by equation~\eqref{e:cocy2} in Theorem~\ref{t:cocy}
we have, for any $(x, 0, y) \in \sG$,
\begin{equation} \label{e:cocy4}
  c(x, 0, y) = f(x) - f(y) + c(\sS{x}, 0, \sS{y}) .
\end{equation}
Let $y$ be chosen to have a dense equivalence class in $\sG(0)$; i.e., so that
$\{ x \mid (x, 0, y) \in \sG(0)\}$ is dense in $P$.  
Replacing $f$ by $g = f + L$ for a constant $L$, we can assume
that $g(y) = 0$. The values of $c_f$ and $c_g$ agree on $\sG(0)$. But
then from equation~\eqref{e:cocy4}
 we have that
\[ 
g(x) = c(x, 0, y) - c(\sS{x}, 0, \sS{y}),
\]
so that $g$ is $\bbZ$-valued on a dense set, hence $\bbZ$-valued on
$P$. 
This completes the first statement.

The proof of the second statement is achieved, 
\emph{mutatis  mutandis}, 
as in the last part of the proof of
Proposition~\ref{p:bdd_cocy}.
\end{proof}

\begin{remark} \label{r:coord}
 Let $G$ be a range finite directed graph with no sources. 
 Any compact, open subset of  $P$  is a finite union of
cylinder sets $Z(\alpha)$. The sets $Z(\alpha)$ can be taken disjoint.
Suppose further that $G$ is a finite graph, so
that $P$ is compact, and let $f \in C(P)$ be a function which
assumes only finitely many values. For each
$t$ in the range of $f$, \ $f^{-1}(t)$ is a compact open subset of
$P$, and hence can be written as a finite,
disjoint union of cylinder sets.  It follows that there is a positive
integer $N$ such that, for all $x \in P$,
the value of $f$ at $x$ depends only on the first $N$ `coordinates' of $x$.
\end{remark}

\begin{theorem} \label{t:integer}
Let $G$ be a finite, transitive, directed graph containing 
at least two distinct simple loops.  Suppose that
$c$ is a $\bbZ$-valued cocycle defined on the AF 
subgroupoid $\sG(0) \subset \sG$, which vanishes precisely on
the unit space $\sG^0$.  Then $c$ has no extension 
to a $\bbZ$-valued cocycle on $\sG$.

If $\sG(0)$ has a dense equivalence class, then $c$ admits no extension to $\sG$.
\end{theorem}

\begin{proof} Let the two simple loops 
be denoted $\alpha$ and $\beta$.
Each of $\alpha$ and $\beta$ contains an edge not
in the other.  By transitivity, there is a path $\gamma$ with 
$r(\gamma) = r(\alpha) = s(\alpha)$ and
$s(\gamma) = r(\beta) = s(\beta)$, and another path $\gamma'$ with 
$r(\gamma') = r(\beta)$, $s(\gamma') = r(\alpha)$. 
Denote the loop $\gamma \beta \gamma'$ by $\beta$. 
$\beta$ may no longer be simple, but it contains an
edge not in $\alpha$, and both loops $\alpha$ and $\beta$ have 
the same initial and terminal vertex. 
Suppose $\alpha$ has $k$ edges, and $\beta$ has $l$ edges.  
If $k \neq l$, we can replace $\alpha$ by
$\alpha^l = \alpha \dots \alpha$ ($l$ concatenations), 
and replace $\beta$ by $\beta^k$.  Changing notation
and denoting the new loops by $\alpha$ and $\beta$, neither 
may be simple, but both now have the same number
of edges, and $\beta$ contains an edge not in $\alpha$.

We suppose $c$ is extendible; hence, there is a function 
$f \in C(P)$ such that the cocycle $c$
is the restriction of $c_f$ to $\sG(0)$. Thus,
\[ 
c(x, 0, y) = \sum_{j=0}^{\infty} [f(\sS^j{x}) - f(\sS^j{y})],
 \text{  for all } (x, 0, y) \in \sG(0). 
\]
Since $c_f$ is assumed to be $\bbZ$-valued, it follows from 
$f(x) = c_f(x, 1, \sS{x})$ that $f$ is $\bbZ$-valued.
As $P$ is compact, $f$ takes on only finitely many values, and
 hence the value of $f$ at $x \in P$ depends
only on some initial path of $x$: $\exists N \in \bbZ^+$ such 
that if $\eta$ is a finite path, $|\eta| \geq N$
and $z, z' \in P$ with $s(\eta) = r(z) = r(z')$, then 
$ f(\eta z) = f(\eta z')$.  (Cf. Remark~\ref{r:coord})

Say $\alpha = e_1\dots e_k$, $\beta = f_1\dots f_k$, $e_i, f_i \in E$. 
By increasing $N$ if necessary, we may assume $k|N$ and that
$\ell = N/k > 1$. 
  Define points
$x, y \in P$ by
\[ 
x = \alpha^{\ell}\alpha \beta \alpha^{\ell} \beta^{\infty} \text{ and }
   y = \alpha^{\ell} \beta \alpha \alpha^{\ell} \beta^{\infty} .
\]
For $z \in P$, let $(z)_N$ denote the truncation of $z$: $(z)_N = z_1 \dots z_N$.
For $n \in \bbZ^+$,  let $m, p$ be determined by the Euclidean
algorithm: $n = mk + p$ with $0 \leq p < k$. 

 Observe that if $m= 0$ (so $n = p$), then
\begin{align*}
(\sS^n{x})_N &= e_{p+1}\dots e_k \alpha^{\ell -1}e_1\dots e_p, \\
(\sS^n{y})_N &= e_{p+1}\dots e_k \alpha^{\ell -1}f_1\dots f_p. 
\end{align*}
If $m = 1$ (so $n = k + p$),
\begin{align*}
(\sS^n{x})_N &= e_{p+1}\dots e_k \alpha^{\ell -1}f_1\dots f_p, \\
(\sS^n{y})_N &= e_{p+1}\dots e_k \alpha^{\ell -2}\beta e_1\dots e_p.
\end{align*} 
For any $1 < m \leq \ell - 1$, with $\alpha^0$  understood to be the empty string,
\begin{align*}
(\sS^n{x})_N &= e_{p+1}\dots e_k \alpha^{\ell -m-1}\alpha \beta \alpha^{m-2}e_1\dots e_p, \\
(\sS^n{y})_N &= e_{p+1}\dots e_k \alpha^{\ell -m-1}\beta \alpha \alpha^{m-2}e_1\dots e_p.
\end{align*}
For $m = \ell$,
\begin{align*}
(\sS^n{x})_N &= e_{p+1}\dots e_k \beta \alpha^{\ell -2}e_1\dots e_p, \\
(\sS^n{y})_N &= f_{p+1}\dots f_k \alpha^{\ell -1}e_1\dots e_p.
\end{align*}
For $ m = \ell + 1$,
\begin{align*}
(\sS^n{x})_N &= f_{p+1}\dots f_k \alpha^{\ell -1}e_1\dots e_p, \\
(\sS^n{y})_N &= e_{p+1}\dots e_k \alpha^{\ell -1}e_1\dots e_p.
\end{align*}
Note that for $m \geq \ell + 2,\ \sS^n{x} = \sS^n{y}$.  Also, 
viewing $n$ as a function of $m$ with $p$
constant, observe that $(\sS^{n(m)}{y})_N = (\sS^{n(m+1)}{x})_N$ 
for $ 0 \leq m  \leq \ell$ and that
$(\sS^{n(0)}{x})_N = (\sS^{n(\ell +1)}y)_N$.  It follows that
\[ 
\sum_{n=0}^{\infty} [f(\sS^n{x}) - f(\sS^n{y})] = 0.
\]
In other words, $c_f(x, 0, y) = 0$; this is impossible since, with 
$x \neq y$, $(x, 0, y)$ is not a unit and
$c_f(x, 0, y) = c(x, 0, y)$ was assumed to vanish only on the unit space.

If now $\sG(0)$ contains a dense orbit, then the second statement of 
the theorem follows immediately from the
second statement of Lemma~\ref{l:integer}
\end{proof}

\begin{remark}
Lemma~\ref{l:integer} and Theorem~\ref{t:integer} together show that
the standard cocycle on $\UHF(n^{\infty})$ (viewed as the core
AF-subalgebra of the Cuntz algebra $O_n$) has no extension to a
cocycle on $O_n$.
\end{remark}

\subsection{The `Analytic' Subalgebra Associated with a Cocycle}
Let $c \in Z^1(\sG, \bbR)$, and set
\[
 A(c) = \{ f \in  C^*(\sG) \mid f \text{ is supported on the
   set } c^{-1}([0,\infty)) \} 
\] 
From Remark~\ref{r:Hinfty}, $f \in A(c)$ if, and only if, the 
maps $t \to \rho(\eta_t(f))$ is an $H^{\infty}$ function
(for $\rho$, $\eta$ as in Remark~\ref{r:Hinfty}).  
Thus $A(c)$ is also written as $H^{\infty}(c)$. 

Clearly, any point in the interior of $c^{-1}([0, \infty))$ 
lies in the spectrum of $A(c)$. On the other hand, 
since the spectrum of any bimodule is open, it follows that 
$\sigma(A(c))$ is the interior of $c^{-1}([0, \infty))$.

Of particular interest are the cocycles $c \in Z^1_0$, that is, 
those which vanish precisely on the unit space, for in that
case we have $\sigma(A(c)) = c^{-1}([0, \infty))$.  If, furthermore, 
the directed graph $G$ satisfies the condition
that every loop has an entrance, then C$^*(\sG^0)$ is a masa in
C$^*(\sG)$ and $A(c)$ is triangular (since
$A(c) \cap A(c)^* = \text{C}^*(\sG^0)$). Furthermore,
$c^{-1}([0, \infty))$ clopen also implies that $A(c) + A(c)^*$ is
dense in C$^*(\sG)$.  Indeed, 
if $\chi_1$ is the characteristic function of $c^{-1}([0, \infty))$
and  $\chi_2$ is the characteristic function of
$c^{-1}((-\infty, 0))$ then any
$f \in C_c(\sG)$ can be written as $f =  f\chi_1 + f\chi_2$.

\section{The Spectral Theorem for Bimodules -- Part III}
\label{s:stb3}  

In this section we extend the spectral theorem for bimodules to show 
that the condition of invariance under
the gauge  automorphisms can be replaced by invariance under the
automorphism group associated with
an `arbitrary' locally constant cocycle (satisfying a mild
constraint). 
As usual, we assume  that the graph $G$ is range
finite and has no sources.

\begin{remark} \label{r:equivnorms}
If $K$ is any compact subset of the groupoid $\sG$, then as Banach
spaces $C(K) \subset \sA$. Since $C(K)$ is complete in both 
the \cstar and supremum norms, these
norms are equivalent on $C(K)$.
\end{remark}

Let $Z(\alpha, \beta)$ be a basic open set in $\sG$. 
 We define a partial homeomorphism $\tau$ on $P$
with $\dom(\tau) = \{\alpha z \in P \mid r(z) = s(\alpha) =
s(\beta) \}$ 
by $\tau(\alpha z) = \beta z $.
By definition, $Z(\alpha, \beta) \subset \sG(k)$, 
where $k = |\alpha| - |\beta|$.

\begin{notation} \label{n:graph} 
With $\tau$ as above, denote
\[ \sG(\tau) = \{ (x, \ell, y) \in \sG:\ y = \tau(x) \} . \]
We refer to $\sG(\tau)$ as the $\sG$-\emph{graph} of $\tau$.
\end{notation}

Note that the graph of $\tau$ could contain points
$(x, \ell, y)$ with $\ell \neq k$.  For example, if $\alpha$ 
is a loop and  $\ell = |\alpha^2| - |\beta|$,
then 
$(\alpha^{\infty}, \ell, \beta \alpha^{\infty})$ 
also lies in $\sG(\tau)$.

\begin{notation}
For $f \in C^*(\sG)$, we let $f_{\tau}$ denote the restriction 
of $f$ to the $\sG$-graph of the partial
homeomorphism $\tau$.  Viewing $f$ as a function on the groupoid 
$\sG$, the restriction is well-defined as a
function on $\sG$.
\end{notation}

Given $f \in \sA$, it is not clear that the restriction 
$f_{\tau}$ also belongs to $\sA$, 
much less that if $f$ belongs to a norm-closed $\sD$-bimodule $\sB$, then 
$f_{\tau} $ also belongs to $\sB$. Our first goal is
to verify these statements.

By~\cite[Cor.~5.5]{MR98g:46083}, path space $P$ is metrizable.  Fix a metric on $P$.
With $\tau$ as above, the domain of $\tau$ is the open 
compact neighborhood $Z(\alpha) \subset P$.
For each $x \in \dom(\tau)$, let $U_n(x)$ be a clopen 
neighborhood centered at $x$ with radius at most
$1/n$.  By compactness, there is a finite subcover, 
$U_n(x_1), \dots, U_n(x_{r_n})$ of dom$(\tau)$. Let 
$U_{n, 1} = U_n(x_1)$ and 
$\ds U_{n, j} = U_n(x_j) \setminus \cup_{i=1}^{j-1} U_n(x_i)$ 
for $j = 2, \dots, r_n$.  Thus the sets
$U_{n, j}$ form a disjoint clopen cover of $Z(\alpha) = \dom(\tau)$.
Let $\chi_{n, j}$ denote the characteristic function of $U_{n,j}$.
Define $\Psi_n \colon \sA \to \sA$ by
\[ 
 \Psi_n(f) =
\sum_{j=1}^{r_n} \chi_{n,j}\cdot f \cdot (\chi_{n,j}\circ \tau^{-1}).
\]
(We identify any function $g \in C_0(P)$ with a function, also denoted
by $g$, on $\sG$ by using the natural identification of $P$ with the
unit space of $\sG$.  The extended function $g$ vanishes at any point
$(x,k,y)$ for which $k \neq 0$ or $x \neq y$.  This function is in
$\sA$; in fact, it is in $\sD$.)

We will show that the sequence $\{ \Psi_n(f)\}_{n=1}^{\infty}$
 converges to $f_{\tau}$, for any
$f \in \sA$. If $f$ happens to be supported on some
 compact subset $K \subset \sG$, then
$\Psi_n(f)$ is also supported on $K$.  Furthermore, if
$(x, k, \tau(x)) \in K\cap \sG(\tau))$, then
$\Psi_n(f)(x, k, \tau(x)) = f(x, k, \tau(x))$ while for 
$ (x, k, y) \notin \sG(\tau)$ we have $\Psi_n(f)(x, k, y) \to 0$.  
The convergence is uniform on compact subsets 
and hence uniform on $\sG$.  By Remark~\ref{r:equivnorms},
$ \Psi_n(f) \to f_{\tau}$  in the \cstar norm.
This, of course, shows tht $f_{\tau} \in \sA$, at least when
$f$ is compactly supported.

To handle the general case, we need to observe that each $\Psi_n$ is norm
decreasing.  Indeed, $\Psi_n$ has the form
$f \mapsto \sum p_n f q_n$ where the sum is finite and each of 
$\{p_n\}$ and $\{q_n\}$ is a family of mutually orthogonal
projections.  Maps of this form on a \cstar algebra are always
contractive. 

\begin{lemma} \label{l:psi}
For $f\in \sA$, $\Psi_n(f)$ converges to $f_{\tau}$. In particular,
$f_{\tau} \in \sA$ and $||f_{\tau}|| \leq ||f||$. If $f$ 
belongs to a norm closed $\sD$-bimodule $\sB$,
then so does $f_{\tau}$.
\end{lemma}

\begin{proof}
We already know that
 $\Psi_n(f)$ converges 
pointwise to $f_{\tau}$ on $\sG$.
Given $\epsilon > 0$, let $g \in \sA$ have compact support, with
$||f - g|| < \epsilon$.  Since $g$ is compactly supported,
 $\Psi_n(g) \to g_{\tau}$ in \cstar norm.  Hence, there is a positive
integer $N$ such that 
 $ ||\Psi_n(g) - g_{\tau}|| < \epsilon$ for all  $n \geq N$.
Then, since $\Psi_n$ is contractive, 
\[
 ||\Psi_n(f) - \Psi_m(f)|| \leq ||\Psi_n(f - g)|| + 
||\Psi_n(g) - \Psi_m(g)|| + ||\Psi_m(g - f)||  < 4\epsilon
\]
for $n, m \geq N$.
  Thus $\{ \Psi_n(f)\}_{n=1}^{\infty}$ 
has a limit, which must agree with its
pointwise limit, $f_{\tau}$. Hence $f_{\tau} \in \sA$ and
\[
||f_{\tau}|| \leq \lim_{n \to \infty} ||\Psi_n(f)|| \leq ||f||.
\]

Finally, it is clear that if $f$ belongs to a 
norm closed $\sD$-bimodule $\sB$, then 
so does each $\Psi_n(f)$; therefore
$f_{\tau} \in \sB$.
\end{proof}

In Theorems~\ref{t:stb_prelim} and~\ref{STB} we proved that a closed
$\sD$-bimodule is determined by its spectrum if, and only if, it is
invariant under the gauge automorphisms.  As noted in 
section~\ref{s:cocycles}, the gauge automorphisms arise in a natural
way from the cocycle $c(x,k,y)=k$ on $\sG$.  In Theorem~\ref{t:stbext}
below we show that a closed $\sD$-bimodule is determined by its
spectrum if, and only if, it is invariant under that one parameter
automorphism group associated with a locally constant cocycle $c$ for
which $c^{-1}(0) \subseteq \sG(0)$.  For a locally constant cocycle
$c$, $c^{-1}(0) \subseteq \sG(0)$ if, and only if, the fixed point
algebra for this one parameter automorphism group is contained 
in the core AF algebra.

As we saw in Theorem~\ref{t:cocy}, continuous cocycles are in
one-to-one correspondence with continuous functions on path space.
Proposition~\ref{p:loc_const} showed that locally constant cocycles
arise from locally constant functions on path space.

Suppose that $\sB$ is a norm closed $\sD$-bimodule.  It is automatic
that $\sB \subseteq A(\sigma(\sB))$; if every Cuntz-Krieger partial
isometry $\SSadj {\alpha}{\beta}$ in $A(\sigma(\sB))$ lies in $\sB$,
then $A(\sigma(\sB)) \subseteq \sB$ and $\sB$ is determined by its
spectrum.  Since we are using the groupoid model, 
$\SSadj {\alpha}{\beta}$ is the characteristic function of
the basic open subset $Z(\alpha,\beta)$ of $\sG$.

A simple observation is useful in the proof of 
Theorem~\ref{t:stbext} below.  Suppose that for each
$(x,k,y) \in \sigma(\sB)$ there is a basic neighborhood 
$Z(\alpha,\beta)$ of $(x,k,y)$ such that
$\SSadj {\alpha}{\beta} \in \sB$.  Then it follows that every
Cuntz-Krieger partial isometry in $A(\sigma(\sB))$ is in $\sB$.
Indeed, if $\SSadj {\alpha}{\beta} \in \sB$ and
$Z(\gamma,\delta) \subseteq Z(\alpha,\beta)$, then
$\SSadj {\gamma}{\delta}$ can be obtained from
$\SSadj {\alpha}{\beta}$ by left and right multiplication by
projections in $\sD$ (use the range projections for
$S_{\gamma}$ and $S_{\delta}$); therefore
$\SSadj {\gamma}{\delta} \in \sB$ also.
If $Z(\nu,\mu)$ is an arbitrary basic open subset of $\sigma(\sB)$,
then by hypothesis, it can be covered by sets $Z(\alpha,\beta)$ for
which $\SSadj {\alpha}{\beta} \in \sB$.  Since $Z(\nu,\mu)$ is
compact, there is a finite subcover.  The observation about subsets
allows us to find a finite subcover of disjoint sets of the form
$Z(\gamma,\delta)$ with $\SSadj {\gamma}{\delta} \in \sB$.  It now
follows that $\SSadj {\nu}{\mu}$ is a finite sum of elements of $\sB$
and so is in $\sB$ itself.

Let $c$ be a real valued cocycle on $\sG$.  Recall that the associated
one parameter automorphism group on $\sA$ is defined by
\[
\eta_z(f)(x,k,y) = z^{c(x,k,y)}f(x,k,y), \quad \text{all } z \in
\bbT. 
\]
(To avoid ambiguity, when $z = e^{it}$ with $0 \leq t < 2\pi$ and $a$
is a real number, we take
$z^a = e^{ita}$.)

\begin{theorem} \label{t:stbext}
Let $\sB \subseteq \sA$ be a norm closed $\sD$-bimodule and let $c$ be a
locally constant cocycle on $\sG$ such that 
the fixed point algebra of the associated one parameter automorphism
group $\eta$ is contained in the core AF algebra.
 Then $\sB = A(\sigma(\sB)$ if, and only
if,  $\sB$ is invariant under $\eta$.
\end{theorem}

\begin{proof}
If $f \in A(\sigma(\sB))$ then $f$ is supported on $\sigma(\sB)$;
clearly each $\eta_z(f)$ is also supported on $\sigma(\sB)$.  Thus, 
$\sB = A(\sigma(\sB)$ trivially implies that $\sB$ is invariant under
the $\eta_z$.

Now assume that $\sB$ is invariant under the $\eta_z$.  By the
observations preceding the theorem, it suffices to prove that for each
point $(x,k,y) \in \sigma(\sB)$, there is a basic open neighborhood
$Z(\alpha,\beta)$ of $(x,k,y)$ for which 
$\SSadj {\alpha}{\beta} \in \sB$.

Given $(x_0,k_0,y_0) \in \sigma(\sB)$ there is an element $f \in \sB$
and a basic neighborhood $Z(\alpha, \beta)$ such that 
$c$ is constant on $Z(\alpha,\beta)$ and
$f(x,k,y) \neq 0$ for all $(x,k,y) \in Z(\alpha,\beta)$.
We will show that $\SSadj {\alpha}{\beta} \in \sB$.

Let $\tau$ be the partial homeomorphism on $P$ with
$\dom(\tau) = Z(\alpha)$, $\ran(\tau) = Z(\beta)$ given by 
$\tau(\alpha z) = \beta z$, for all $\alpha z \in Z(\alpha)$.
By Lemma~\ref{l:psi}, $f_{\tau} \in \sB$.

Let $a=c(x_0, k_0, y_0)$.  Define $E \colon \sA \to \sA$ by
\[
E(g)(x,k,y) = \int_{\bbT} z^{-a} \eta_z(g) (x,k,y)\,dz.
\]
The integration is with respect to normalized Lebesgue measure on
$\bbT$. Each $\eta_z$ is isometric, so $E$ is contractive.

Since $f_{\tau} \in \sB$, $E(f_{\tau})\in \sB$ also.  When
$(x,k_0,y) \in Z(\alpha,\beta)$, we have
$z^{-a}\eta_z(f_{\tau})(x,k_0,y) = f_{\tau}(x,k_0,y)$.  If
$(x,k_0,y) \notin Z(\alpha,\beta)$, then
$z^{-a}\eta_z(f_{\tau})(x,k_0,y) = f_{\tau}(x,k_0,y)$ again holds,
since both sides of the equation equal 0.  Thus,
$E(f_{\tau})(x,k_0,y) = f_{\tau}(x,k_0,y)$, for all
$(x,k_0,y) \in \sG$.  When $k \neq k_0$, then
\[
c(x,k,y) -a = c(x,k,y) - c(x,k_0,y) = c(x,k,y) +c(y,-k_0,x) = 
c(x,k-k_0,y) \neq 0.
\]
(The inequality follows from the fact that
$c^{-1}(0) \subseteq \sG(0)$.)
Therefore $a \neq c(x,k,y)$ and the integrand in
 $E(f_{\tau})(x,k,y)$ is a
non-zero power of $z$ multiplied by $f_{\tau}(x,k,y)$.  It follows
that $|E(f_{\tau})(x,k,y)| < |f_{\tau}(x,k,y)|$ with the ratio between
the two numbers dependent only on $(x,k,y)$.

If we now let $E^n$ be the $n$-fold composition of $E$ with itself, we
have $E^n(f_{\tau})(x,k_0,y) = f_{\tau}(x,k_0,y)$ for all
$(x,k_0,y) \in Z(\alpha,\beta)$ and
$E^n(f_{\tau})(x,k,y) \to 0$ otherwise.

Let $f_{Z(\alpha,\beta)}$ denote the restriction of $f_{\tau}$ to 
$Z(\alpha,\beta)$.  If $f_{\tau}$ has compact support, then
$E^n(f_{\tau}) \to f_{Z(\alpha,\beta)}$ uniformly and (by
Remark~\ref{r:equivnorms}) in \cstar norm as well.  It follows in this
case that $f_{Z(\alpha,\beta)} \in \sB$.

For the general case, let $\epsilon > 0$ and let $g \in \sA$ have
compact support and satisfy $\|f - g\| < \epsilon$.  Then
$g_{\tau}$ has compact support and
$\|f_{\tau} - g_{\tau}\| \leq \|f -g\| < \epsilon$ 
(Lemma~\ref{l:psi}).  With $g_{Z(\alpha, \beta)}$ the restriction of
$g$ to $Z(\alpha,\beta)$, we know that there is $N \in \bbN$ such that
$\|E^n (g_{\tau} - g_{Z(\alpha, \beta)})\| < \epsilon$ for all 
$n \geq N$.  Therefore, when $n,m \geq N$,
\[
\|E^n(f_{\tau}) - E^m(f_{\tau})\| \leq 
\| E^n(f_{\tau}- g_{\tau})\| + 
\|E^n(g_{\tau}) - E^m(g_{\tau})\| + \| E^m(g_{\tau}- f_{\tau})\| 
< 4\epsilon.
\]
Thus, even when $f_{\tau}$ is not compactly supported,
$E^n(f_{\tau})$ is convergent in \cstar norm; the limit must agree
with the pointwise limit $f_{Z(\alpha,\beta)}$.  Since
each $E^n(f_{\tau}) \in \sB$, we obtain 
$f_{Z(\alpha,\beta)} \in \sB$.

Now define a continuous, compactly supported function $h$ on $P$ by
\[
h(x) =
\begin{cases}
\frac 1 {f(x,k_0,y)}, \quad &x \in \dom (\tau), \\
0, &\text{otherwise}.
\end{cases}
\]
Extending $h$ to all of $\sG$ by taking it to have value 0 off $P$, we
view $h$ as an element of $\sD$.  But now
$\SSadj {\alpha}{\beta} = h f_{Z(\alpha,\beta)} \in \sB$.  This
completes the proof.
\end{proof}

\section{Nest Subalgebras of Graph \cstar algebras} \label{s:nsa}

An additional structure on on $G$ -- a total ordering of the edges --
leads in a natural way to nest subalgebras of $\sA$.
Arbitrary total orders on $E$ appear to be too general, so we consider
orders on $E$ which are compatible with a total ordering on $V$.
Given a total order on $V$ and, for each $v \in V$,  a total order on 
$\{e \mid r(e) = v\}$, we can then define an order on $E$ in which two
edges with the same range are ordered by the order on 
$\{e \mid r(e) = v \}$ and two edges with different ranges are ordered
by the order on $V$.  We could, of course, use the sources
instead of the ranges, or even combine the two; but it is orders
compatible with the ranges which are most suitable for the algebras
which we shall study.  There is a way to rephrase the definition of
the orders we study; we use this for the formal definition:

\begin{definition}
An \emph{ordered graph} is a directed graph $G$ together with a total
order $\preceq$ on $E$ which satisfies the property that, for each 
$v \in V$, $\{e \mid r(e) =v\}$ is an interval in the order on $E$.
\end{definition}

Throughout this section we assume that
 $G$ is a finite  ordered graph.
  We use the order on the graph to define a nest of
 projections in $\sD$; the associated nest subalgebra of $\sA$ is the
 object of study.

For each $k$, the (left to right) lexicographic order 
gives a total order on $F_k$, the set of paths of length $k$.  
(The lexicographic order is based on the order on $E$.)
We denote this order by $\preceq$.
For each finite path $\alpha = \alpha_1 \dots \alpha_k$,
let $R_{\alpha}$ denote the range of the partial isometry
$S_{\alpha} = S_{\alpha_1} \dots S_{\alpha_k}$.  
$\{R_{\alpha} \mid |\alpha| = k\}$ is an orthogonal set of
projections which sum to the identity.  This set inherits a total
order from the  lexicographic order on 
$F_k$.  We shall use the notation $\ll$ for the strict variant of
this total order.

Let $\sN_k$ be the nest whose atoms are the $R_{\alpha}$ with
$|\alpha|=k$, taken in the order above.  
Projections in $\sN_k$ have the form 
$\sum R_{\alpha}$, summed over initial segments in the order $\ll$.

Let $\alpha = \alpha_1 \dots \alpha_k \in F_k$.   Write
$\{e \in E \mid r(e) = s(\alpha_k)\}$ as $\{e_1, \dots, e_p \}$ with
$e_1 \prec e_2 \prec \dots \prec e_p$.  Then
$\{\alpha e_1, \dots, \alpha e_p \}$ forms an order interval
in $F_{k+1}$.
For each path $\alpha$, $R_{\alpha} =\sum R_{\alpha e} $, where the sum is
over all edges $e$ such that $r(e) = s(\alpha)$. 
It follows that
$\sN_k \subseteq \sN_{k+1}$, for all $k$.
Let $\sN = \bigcup \sN_k$ and 
$\Alg \sN = \{ A\in \sA \mid AP=PAP \text{ for all }P \in \sN\}$.

Note that, for fixed $k$ and $p$,
 $\sN = \bigcup_n \sN_{p+nk}$.  Consequently,
to prove that an element $A$ of the graph \cstar algebra is in
$\Alg \sN$ it suffices to prove that $A \in \Alg \sN_{p+nk}$, for all
$n$. 

\begin{definition}
We shall refer to $\sN$ and $\Alg \sN$ as the \emph{nest} and the
\emph{nest algebra} induced by the order $\preceq$ on $E$.
\end{definition}

\begin{remark}
The material in this section was inspired by, and is an extension of,
the work on the Volterra subalgebra of the Cuntz \cstar algebra found
in~\cite{MR86d:47057} and~\cite{hop_peters}.  
The Cuntz algebra $O_n$ is the graph \cstar
algebra for a graph with one vertex and $n$ loops.  By symmetry, 
the choice of order on the $n$ loops is immaterial .    (Indeed, it is
not hard to find a unitary in $O_n$ which conjugates
the generators  in one order into
the generators in another order.)

There is a natural representation of $O_n$ acting on $L^2[0,1]$.  For
each $k = 1, \dots, n$, let
$S_k$ be the isometry on
$L^2[0,1]$  associated
 with the affine, order preserving
map from $[0,1]$ onto the interval
$\ds \left[\frac{k-1}n, \frac kn\right]$.
  If the $n$ loops in
the graph for $O_n$ are $e_1, \dots, e_n$ in order, then 
$S_1, \dots, S_n$ are the corresponding generating isometries.  The nest
$\sN$ then consists of the projections which correspond to the
intervals of the form $\ds \left[0, \frac j{n^k}\right]$, 
where $j$ and $k$ are
non-negative integers.  This nest is strongly dense in the Volterra
nest (which consists of projections corresponding to intervals
$[0,t]$, $0 \leq t \leq 1)$.  In this representation, the nest algebra
$\Alg \sN$ is exactly the intersection of $O_n$ with the usual
Volterra nest algebra acting on $L^2[0,1]$.
\end{remark}

Observe that $\Alg \sN$ is invariant under the gauge automorphisms.
Indeed, if $z \in \bbT$, $A \in \Alg \sN$ and $P \in \sN$, then,
since $P$ is in the fixed point algebra of the gauge automorphisms,
$P^{\perp}\eta_z(A)P = \eta_z(P^{\perp}AP) = 0$.  Thus,
$\eta_z(A) \in \Alg \sN$, for all $A \in \Alg \sN$, $z \in \bbT$.

By the spectral theorem for bimodules (Theorem~\ref{STB}), 
$\Alg \sN$ is the closed linear span of the Cuntz-Krieger partial
isometries which it contains.   We will now characterize
the Cuntz-Krieger partial isometries $\SSadj{\alpha}{\beta}$
 in $\Alg \sN$ in terms of the properties of the finite paths
$\alpha$ and $\beta$.  
This, in turn, will enable us to give a description of 
the spectrum $\sigma(\Alg \sN)$.

\begin{definition}
A path $\alpha$ is \emph{s-minimal} if $\alpha \preceq \beta$ whenever
$\beta$ is a path with $|\beta|=|\alpha|$ and 
$r(\beta)=s(\alpha)$.
$\alpha$ is \emph{s-maximal} if $\beta \preceq \alpha$ whenever
$\beta$ is a path with $|\beta|=|\alpha|$ and 
$r(\beta)=s(\alpha)$.
\end{definition}

\begin{remark} \label{r:Cuntz_paths}
In a Cuntz algebra $O_n$, finite paths are essentially finite
sequences from the integers $\{1, \dots, n\}$.  A finite path $\alpha$
is $s$-minimal if $\alpha_i=1$ for all $i$ and $s$-maximal if
$\alpha_i=n$ for all $i$.
\end{remark}

\begin{proposition} \label{alphaonly}
$S_{\alpha} \in \Alg \sN \Longleftrightarrow \alpha 
\text{ is $s$-minimal}$.
\end{proposition}

\begin{proof}
Suppose that $\alpha$ is not $s$-minimal.  Then there is a path
 $\beta$  with
$|\beta| = |\alpha|$, $r(\beta)=s(\alpha)$, and
$\beta \prec \alpha$.  With $k$ the common degree of $\alpha$ and
 $\beta$, $R_{\alpha}$ and $R_{\beta}$ are atoms from $\sN_k$ and
$R_{\beta} \ll R_{\alpha}$.  Now, $S_{\alpha}$ is non-zero on
$R_{\beta}$ (since $r(\beta) = s(\alpha)$) and so
$R_{\alpha} S_{\alpha} R_{\beta} \neq 0$.  But then
$S_{\alpha} \notin \Alg \sN_k$ and hence
$S_{\alpha} \notin \Alg \sN$.

Now suppose that  $\alpha$ is $s$-minimal.  We distinguish two cases.
First assume that $r(\alpha) \neq s(\alpha)$; i.e., $\alpha$ is not a
loop. Then the initial space $Q_{\alpha}$ is a sum of atoms of the
  form $R_{\beta}$, where $|\beta|=|\alpha|$ and
$\alpha \prec \beta$.  If we let $P$ be the smallest projection in
$\sN_k$ such that $R_{\alpha} \leq P$, then each of the $R_{\beta}$ in
the sum for $Q_{\alpha}$ is orthogonal to $P$.  Therefore,
$S_{\alpha} = P S_{\alpha} P^{\perp}$ and
$S_{\alpha} \in \sN$.

Next assume that $\alpha$ is a loop.  Then the initial space
$Q_{\alpha}$ can be written as a sum
$R_{\alpha} + \sum R_{\beta}$, where the $\beta$ in the sum run over
 paths with $|\beta|=|\alpha|$, $r(\beta)=s(\alpha)$, and
$\alpha \prec \beta$.  Since
$S_{\alpha} R_{\beta} = R_{\alpha}S_{\alpha} R_{\beta}$ for each such 
$\beta$,  each $S_{\alpha} R_{\beta} \in \Alg \sN$.

It remains to show that $S_{\alpha} R_{\alpha} \in \Alg \sN$. 
Let $k$ be the degree of $\alpha$.   As noted
above, it is sufficient to show that 
$S_{\alpha} R_{\alpha} \in \Alg \sN_{nk}$ for each 
positive integer $n$.
Let $P$ be a projection in $\sN_{nk}$.  If $P \perp R_{\alpha}$, then
$S_{\alpha} R_{\alpha} P=0$ and $S_{\alpha} R_{\alpha}$ 
trivially leaves $P$ invariant.
$P$ is also trivially left invariant if $R_{\alpha} \leq P$.
This leaves the case in which 
$0 < PR_{\alpha} < R_{\alpha}$.  To show that $P$ is invariant under
$S_{\alpha} R_{\alpha}$, it suffices to prove
$S_{\alpha} R_{\alpha} P \in \Alg \sN_{nk}$.

The projection $PR_{\alpha}$ can be written as a sum of
atoms $R_{\beta}$ (from $\sN_{nk}$) where the $\beta$ run through an
interval in the order on $F_{nk}$.  Let $\beta$ be one of
these paths.  Write $\beta = \beta_1 \dots \beta _n$, where
each $\beta_i$ has length $k$.  We need to show that 
$S_{\alpha} R_{\beta}$ has range contained in atoms whose indices
precede or equal $\beta$.  The range of $S_{\alpha}R_{\beta}$ is
$R_{\alpha\beta}$, which is a subprojection of
$R_{\alpha\beta_1 \dots \beta_{n-1}}$.   The assumption that
$\alpha$ is $s$-minimal implies that $\alpha \preceq \beta_1$.
If $\alpha \prec \beta_1$, then
$R_{\alpha\beta_1 \dots \beta_{n-1}} \ll R_{\beta_1 \dots \beta_n}=R_{\beta}$.
On the other hand, if $\beta_1 = \alpha$, then 
$r(\beta_2)=s(\beta_1)=s(\alpha)=r(\alpha)$ and
$\alpha \preceq \beta_2$. (Again use $\alpha$ is $s$-minimal).

 Once again, if $\alpha \prec \beta_2$ then 
$R_{\alpha\beta_1 \dots \beta_{n-1}} =R_{\alpha \alpha \beta_2 \dots \beta{n-1}}
\ll R_{\alpha \beta_2 \dots \beta_n} = R_{\beta_1 \dots \beta_n}$.
Continuing in this fashion, we see that if any of the $\beta_i$ are
unequal to $\alpha$, then $S_{\alpha} R_{\beta} \in \Alg \sN_{nk}$.
Finally, if all $\beta_i = \alpha$, then $S_{\alpha}$ maps
$R_{\alpha \dots \alpha}$ into itself, so again 
$S_{\alpha} R_{\beta}=
S_{\alpha}R_{\alpha \dots \alpha} \in \Alg \sN_{nk}$.
From this it follows that
$S_{\alpha} R_{\alpha} P \in \Alg \sN_{nk}$ and the Proposition is proven.
\end{proof}

\begin{proposition} \label{alphaequalsbeta}
  Let $\alpha$ and $\beta$ be two paths of equal length with
 $s(\alpha) = s(\beta)$ 
\/{\rm (}so that $\SSadj{\alpha}{\beta} \neq 0$\/{\rm )}.  Then
$\SSadj{\alpha}{\beta} \in \Alg \sN$ if, and only if, $\alpha \preceq \beta$.
\end{proposition}

\begin{proof}
If $\alpha = \beta$, then $\SSadj{\alpha}{\beta}$ is a projection in the canonical
diagonal of the graph \cstar algebra.  Since the nest $\sN$ is also in
this diagonal, $\SSadj{\alpha}{\beta} \in \Alg \sN$. Note that
$\SSadj{\alpha}{\beta} = R_{\alpha} \SSadj{\alpha}{\beta} R_{\beta}$.
 If $\alpha \prec \beta$, then
$R_{\alpha} \ll R_{\beta}$ and $\SSadj{\alpha}{\beta} \in \Alg \sN$.  If
$\beta \prec \alpha$ then $R_{\beta} \ll R_{\alpha}$ and
$\SSadj{\alpha}{\beta} \notin \Alg \sN$.
\end{proof}

\begin{proposition} \label{betashorteralpha}
Let $\alpha$ and $\beta$ be two paths with 
$0<k = |\beta| < |\alpha|$ and $s(\alpha) = s(\beta)$.  Write
$\alpha = \delta \gamma$ where $|\delta| = |\beta|$ and
$|\gamma| = |\alpha| - |\beta|$.  Then
$\SSadj{\alpha}{\beta} \in \Alg \sN$ if, and only if, one of the following two
conditions holds\/{\rm :}
\begin{enumerate}
\item $\delta \prec \beta$,
\item $\delta = \beta$ and $\gamma$ is $s$-minimal.
\end{enumerate}
\end{proposition}

\begin{proof}
 The initial space for
$\SSadj{\alpha}{\beta}$ is the final space for $S_{\beta}$, namely $R_{\beta}$.    The
range space is $R_{\alpha} = R_{\delta \gamma}$, which 
is a subprojection of
$R_{\delta}$. If  $\delta \prec \beta$ then
 $R_{\delta} \ll R_{\beta}$ as
atoms from $\sN_k$ and $\SSadj{\alpha}{\beta} \in \Alg \sN$.  And  if
$\beta \prec \delta$ then $R_{\beta} \ll R_{\delta}$ and
$\SSadj{\alpha}{\beta} \notin \Alg \sN$.

Assume that $\delta = \beta$, so that  $\alpha = \beta \gamma$.
  Observe that $r(\gamma) = s(\beta)$ and that
 $s(\beta) = s(\alpha) = s(\gamma)$.  Thus,
$r(\gamma) = s(\gamma)$ and $\gamma$ is a loop.

Now the initial space for
$S_{\beta}^*$ is $R_{\beta}$ and the final space is
 the initial space for
$S_{\alpha}$.  Therefore, $\SSadj{\alpha}{\beta}$ maps $R_{\beta}$ onto
$R_{\alpha} = R_{\beta \gamma}$, which is a subprojection
of $R_{\beta}$.  So $\SSadj{\alpha}{\beta}$ trivially leaves invariant any projection
which contains $R_{\beta}$ or is orthogonal to $R_{\beta}$.

Assume that $\gamma$ is $s$-minimal.
Let $t= |\gamma|$.  It is sufficient to prove that, for any positive
integer $n$, $\SSadj{\alpha}{\beta} \in \Alg \sN_{k+nt}$.  By the  preceding 
paragraph, it is
enough to look at the action of $\SSadj{\alpha}{\beta}$ on atoms from
 $\sN_{k+nt}$ which
are subprojections of $R_{\beta}$.  Each of these atoms has the form
$R_{\beta \eta_1 \dots \eta_n}$, where $|\eta_i| = t$, for all $i$.
Now $S_{\beta}^*$ maps $R_{\beta \eta_1 \dots \eta_n}$ onto
$R_{\eta_1 \dots \eta_n}$ and so $\SSadj{\alpha}{\beta}$ maps 
$R_{\beta \eta_1 \dots \eta_n}$ onto
$R_{\beta \gamma \eta_1 \dots \eta_{n}}$.  The latter is a
subprojection of
$R_{\beta \gamma \eta_1 \dots \eta_{n-1}}$, which is an atom
from $\sN_{k+nt}$.  So all we have to do is to prove that
$R_{\beta \gamma \eta_1 \dots \eta_{n-1}}$ precedes or equals
$R_{\beta \eta_1 \dots \eta_n}$ in the ordering for atoms from
$\sN_{k+nt}$.

If every $\eta_i = \gamma$, this is trivial.  Otherwise, let
$\eta_j$ be the first $\eta$ which is unequal to $\gamma$. 
If $j=1$, then $r(\eta_1) = s(\gamma)$. If $j>1$, then
$r(\eta_j) = s(\eta_{j-1}) = s(\gamma)$.  Since $\gamma$ is
$s$-minimal, $\gamma \prec \eta_j$.  But then
$R_{\beta \gamma \eta_1 \dots \eta_{n-1}} \ll
R_{\beta \eta_1 \dots \eta_n}$ .

It remains to show that if $\gamma$ is not $s$-minimal then
$\SSadj{\alpha}{\beta} \notin \Alg \sN$.  Suppose that $\eta$ is a path with
$|\eta| = |\gamma|$, $r(\eta) = s(\gamma) = s(\beta)$ and
$\eta \prec \gamma$.  Then $\beta \eta \prec \beta \gamma$ and
$R_{\beta \eta} \ll R_{\beta \gamma}$.  Now $S_{\beta}^*$ maps
$R_{\beta \eta}$ into the initial space for 
$S_{\beta \gamma} = S_{\alpha}$ and $\SSadj{\alpha}{\beta}$ maps 
$R_{\beta \eta}$ into a subprojection of $R_{\beta \gamma}$; thus
$\SSadj{\alpha}{\beta} \notin \sN$.
\end{proof}

\begin{corollary}
Suppose that $\gamma$ is a path such that $S_{\gamma} \in \Alg \sN$.
Then, for any $\beta$, 
$S_{\beta}^{\vphantom{*}} S_{\gamma}^{\vphantom{*}} 
S_{\beta}^* \in \Alg \sN$.
\end{corollary}

\begin{proof}
By Proposition~\ref{alphaonly}, $\gamma$ is $s$-minimal.
If $s(\gamma) \neq s(\beta)$ then
$S_{\beta}^{\vphantom{*}} S_{\gamma}^{\vphantom{*}} 
S_{\beta}^* = 0 \in \Alg \sN $.  Otherwise, condition 2 of
Proposition~\ref{betashorteralpha} yields the corollary.
\end{proof}

The next two propositions can be proven with arguments analogous to
the ones used in Proposition~\ref{alphaonly} and
Proposition~\ref{betashorteralpha}.  However, a shortcut is available.
If we reverse the order on paths of length $k$ and therefore reverse
the order on the corresponding atoms, we obtain the nest
$\sN^{\perp}$ instead.  Since a path is $s$-minimal with respect to the
reversed order if, and only if, it is $s$-maximal with respect to the
original order and since
$\Alg \sN^{\perp} = (\Alg \sN)^*$, Proposition~\ref{betaonly} and
Proposition~\ref{alphashorterbeta} are immediate consequences of
Proposition~\ref{alphaonly} and Proposition~\ref{betashorteralpha}.

\begin{proposition} \label{betaonly}
$S_{\beta}^* \in \Alg \sN \Longleftrightarrow \beta$ is $s$-maximal.
\end{proposition}

\begin{proposition} \label{alphashorterbeta}
Let $\alpha$ and $\beta$ bet two paths with 
$0< |\alpha| < |\beta|$ and $s(\alpha)=s(\beta)$.
Write $\beta = \delta \gamma$ where
$|\delta| = |\alpha|$ and $|\gamma| = |\beta| - |\alpha|$.
Then $\SSadj{\alpha}{\beta} \in \Alg \sN$ if, and only if, one of the following two
conditions holds\/{\rm :}
\begin{enumerate}
\item $\alpha \prec \delta$,
\item $\alpha = \delta$ and $\gamma$ is $s$-maximal.
\end{enumerate}
\end{proposition}

The following theorem summarizes the sequence of Propositions above:

\begin{theorem} \label{CKinNA}
Let $G$ be a finite, ordered graph and let $\sN$ be the associated
nest.  A Cuntz-Krieger partial isometry $\SSadj {\alpha}{\beta}$ lies
in $\Alg \sN$ if, and only if, one of the 
following conditions holds\/{\rm :}
\begin{enumerate}
\item $|\alpha| = |\beta|$ and $\alpha \preceq \beta$. \label{eqlength}
\item $\alpha = \delta \gamma$ with $|\delta| = |\beta|$ and
$\delta \prec \beta$. \label{dprecb}
\item $\alpha = \beta \gamma$ and $\gamma$ is $s$-minimal. \label{gmin}
\item $\beta = \delta \gamma$ with $|\delta| = |\alpha|$ and 
$\alpha \prec \delta$. \label{aprecd}
\item $\beta = \alpha \gamma$ and $\gamma$ is $s$-maximal. \label{gmax}
\end{enumerate}
\end{theorem}

We can now characterize the points $(x,k,y) \in \sG$ 
that are  in the spectrum of $\Alg \sN$.  Note that
path space $P$ is totally ordered by the lexicographic order
based on the total order on $E$; once again we let $\preceq$ denote
this order.  (We will, in fact, need to compare $x$ and $y$ only when
$x$ and $y$ are shift equivalent.)

\begin{theorem} \label{spectrumNA}
Let $G$ be a finite ordered graph and let $\sN$ be the associated nest.
A point $(x,k,y) \in \sG$ is in $\sigma(\Alg \sN)$ if, and only if,
one of the following conditions holds\/{\rm :}
\begin{enumerate}
\item $x \prec y$. \label{radical_pt}
\item  $x = y$ and $k = 0$. \label{diagonal_pt}
\item  $x=y$, $k>0$  and $x=\beta \gamma \gamma \gamma \dots$ where
  $\gamma$ is $s$-minimal and $|\gamma| =k$. \label{min_pt}
\item $x=y$, $k<0$  and $x=\alpha \gamma \gamma \gamma \dots$ where
$\gamma$ is $s$-maximal and $|\gamma| = -k$. \label{max_pt}
\end{enumerate}
\end{theorem}

\begin{proof}
The proof is, of course, based on the fact that
$\sigma(\Alg \sN)$ is the union of the sets
$Z(\alpha,\beta)$ with $\alpha$ and $\beta$ satisfying any of the five
conditions in the Theorem~\ref{CKinNA}. 

First suppose that $(x,k,y) \in \sG$ and
$x \prec y$ in the lexicographic order.
We can then find
$\alpha$ and $\beta$ satisfying one of 
conditions~\ref{eqlength}, \ref{dprecb}, or~\ref{aprecd}
(depending on the value of $k$) in Theorem~\ref{CKinNA}, so that 
$(x,k,y) \in Z(\alpha,\beta) \subseteq \sigma(\Alg \sN)$.
Thus
$\{(x,k,y) \in \sG \mid x \prec y \} \subseteq \sigma(\Alg \sN)$.
Equally well, if $\alpha$ and $\beta$ 
satisfy~\ref{eqlength} (with $\alpha \neq \beta$), 
\ref{dprecb}, or~\ref{aprecd}, and
$(x,k,y) \in Z(\alpha,\beta)$, then $x \prec y$.  

By condition~\ref{eqlength} in Theorem~\ref{CKinNA},
any set of the 
form $Z(\alpha,\alpha)$ is contained in
$\sigma(\Alg \sN)$; thus $(x,0,x) \in \sigma(\Alg \sN)$, for all $x$.

Now suppose that  $(x,k,y) \in Z(\alpha,\beta)$ when
$\alpha = \beta \gamma$ and $\gamma$ is $s$-minimal. (So $k>0$.) 
 Then 
\begin{align*}
  x & = \beta_1 \dots \beta_n \gamma_1 \dots \gamma_k z_1 z_2 \dots 
\text{ and}\\
y &= \beta_1 \dots \beta_n  z_1 z_2 \dots .
\end{align*}
Since $\gamma$ is $s$-minimal and $z_1 \dots z_k$ is a finite path
whose range is the source of $\gamma$, we have
$\gamma \preceq z_1 \dots z_k$.  If $\gamma \prec z_1 \dots z_k$, then
$x \prec y$.  So, suppose that $\gamma = z_1 \dots z_k$.  Now
$z_{k+1} \dots z_{2k}$ is a finite path whose range is the source of
$\gamma$ and so $\gamma \preceq z_{k+1} \dots z_{2k}$.  Again, if
$\prec$ holds, then $x \prec y$; otherwise 
$z_{k+1} \dots z_{2k} = \gamma$. It is now clear that an induction
argument shows that either $x \prec y$ or $x=y$ has the form
$\beta \gamma \gamma \gamma \dots$, where $\gamma$ is
$s$-minimal.  The points $(x,k,y)$ with $x \prec y$ have already been
covered by the previous discussion, so the new points in 
$\sigma(\Alg \sN)$ are the ones of the form $(x,k,x)$ where
$x = \beta \gamma \gamma \gamma \dots$, $k=|\gamma|$, and
$\gamma$ is $s$-minimal.  Any point of $\sG$ of this form is in a
suitable $Z(\alpha, \beta)$ and so is in $\sigma(\Alg \sN)$.

We can analyze $Z(\alpha,\beta)$ when $\beta = \alpha \gamma$ and
$\gamma$ is $s$-maximal in a similar way.  If 
$(x,k,y) \in Z(\alpha, \beta)$ and $x \ne y$ then $x \prec y$.  If
$x=y$ then $-k = |\beta| - |\alpha| = |\gamma|$ and
$x = \alpha \gamma \gamma \gamma \dots$ with $\gamma$ $s$-maximal.
Any point with this form is in $\sigma(\Alg \sN)$.

All that remains is the trivial observation that
if $y \prec x$ then
$(x,k,y) \notin \sigma(\Alg \sN)$.
\end{proof}

We next determine the spectrum of the (Jacobson) radical of a nest
subalgebra of a graph \cstar algebra.  Since the radical is invariant
under automorphisms, we know that it is determined by its spectrum.
Analogy with the case of upper triangular matrices and with refinement
subalgebras of AF \cstar algebras suggests that 
 the spectrum of the radical consists of those points $(x,k,y)$ in
$\sigma(\Alg \sN)$ with $x \prec y$ (condition~\ref{radical_pt} in 
Theorem~\ref{spectrumNA}).  Indeed,

\begin{proposition} \label{spec_rad}
The set  $R = \{(x,k,y) \in \sigma(\Alg \sN) \mid x \prec y \}$
is the spectrum of the radical of $\Alg \sN$.
 Consequently,
$A(R)$ is the radical of $\Alg \sN$.
\end{proposition}

\begin{proof}
Theorem~\ref{STB} implies that the second statement follows from the
  first.  Temporarily, let
$R_0$ denote the spectrum of the radical of $\Alg \sN$.
We need to prove that $R_0 = R$.

We first show that $R \subseteq R_0$.  Let $(x,k,y)\in R$, so that
$x \prec y$.  Choose finite strings $\alpha$ and $\beta$ such that
$s(\alpha) = s(\beta)$, $\alpha \prec \beta$ and
$(x,k,y) \in \sigma(\SSadj {\alpha} {\beta})$.  ($\alpha$ and $\beta$
need not have the same length; by $\alpha \prec \beta$ we simply mean
that there is an index $j$ such that $\alpha_i = \beta_i$ for $i<j$
and $\alpha_j \prec \beta_j$.)  Now, the range projection for 
$\SSadj {\alpha} {\beta}$ is contained in the atom
$R_{\alpha_1 \dots \alpha_j}$ and the initial projection is
contained in $R_{\beta_1 \dots \beta_j}$.  Since
$R_{\alpha_1 \dots \alpha_j} \ll R_{\beta_1 \dots \beta_j}$, it
follows that if $P$ is the smallest projection in $\sN$ which contains
$R_{\alpha_1 \dots \alpha_j}$, then
 $R_{\beta_1 \dots  \beta_j}$ is orthogonal to $P$.  Thus,
$\SSadj {\alpha} {\beta} = P \SSadj {\alpha} {\beta}
P^{\perp}$ and $\SSadj {\alpha} {\beta}$ lies in the radical
of $\Alg \sN$.  Therefore $(x,k,y) \in R_0$ and
$R \subseteq R_0$.

To complete the proof, we need to show that any point of 
$\sigma(\Alg \sN)$ which satisfies conditions~\ref{diagonal_pt},
\ref{min_pt} or~\ref{max_pt} is not in $R_0$.  For points of the form
$(x,0,x)$ this is trivial -- they are in the suppport set of a
non-zero projection and the radical contains no non-zero projections.
The arguments for points which satisfy conditions~\ref{min_pt} and
\ref{max_pt} are very similar, so we just treat the first of these two
cases.

Assume $k>0$, $x = \delta \gamma \gamma \gamma \dots$, $k = |\gamma|$,
and $\gamma$ is $s$-minimal.  Suppose that
$(x,k,x) \in R_0$.  Since $R_0$ is open, there is a positive integer
$n$ so that if
\begin{alignat*}{2}
\alpha &=  \delta \gamma \dots \gamma   &(n+1 \text{ copies of }
\gamma) \\
\beta &= \delta \gamma \dots \gamma \qquad &(n \text{ copies of } \gamma)
\end{alignat*}
then $Z(\alpha, \beta) \subseteq R_0$.  Since $r(\gamma) = s(\delta)$
and $r(\gamma) = s(\gamma)$, it follows from the Cuntz-Krieger
relations that $\ssadj{\gamma} \leq \sadjs{\delta}$ and
$\ssadj{\gamma} \leq \sadjs{\gamma}$.  Using this, and the fact that
$\alpha$  has one more copy of $\gamma$ than $\beta$ has, we obtain
\begin{displaymath}
  \SadjS{\beta}{\alpha} = S_{\gamma}^* \dots S_{\gamma}^* S_{\delta}^*
S_{\delta}^{\vphantom{*}} S_{\gamma}^{\vphantom{*}} \dots
S_{\gamma}^{\vphantom{*}}  = S_{\gamma}.
\end{displaymath}
and, hence,
  $(\SSadj{\alpha}{\beta})^2 = S_{\alpha}^{\vphantom{*}}
S_{\gamma}^{\vphantom{*}}S_{\beta}^*$.  Since 
$s(\beta) = s(\gamma)$ and $r(\gamma)=s(\alpha)$, 
$(\SSadj{\alpha}{\beta})^2 \ne 0$.  Note that 
$\alpha \gamma$ has the same form as $\alpha$ except that there are
now $n+2$ $\gamma$'s feeding into  $\delta$.  

The same considerations as above show that if
$\alpha^{(p)} = \delta \gamma \dots \gamma$ with $n+p$ copies of
$\gamma$, then 
$\SadjS{\beta}{{\alpha^{(p)}}} = S_{\gamma}\dots S_{\gamma}$ 
($p$ copies of $\gamma$) and the square of
$(\SSadj {{\alpha^{(p)}}}{\beta})^2 $
is a non-zero partial isometry and so  has norm 1.
At this point it is now a simple matter to show that
$\|(\SSadj{\alpha}{\beta})^k\| = 1$ for all $k$ and therefore
that $\SSadj{\alpha}{\beta}$ is not quasi-nilpotent.  But then
$\SSadj{\alpha}{\beta}$ is not in the radical of
$\Alg \sN$, contradicting $Z(\alpha,\beta) \subseteq R_0$.
Thus any point in the spectrum of $\Alg \sN$ which satisfies condition
\ref{min_pt} of Theorem~\ref{spec_rad} is not in $R_0$.  As mentioned
earlier, points sstisfying condition~\ref{max_pt} are handled similarly.
\end{proof}

It was shown in~\cite{MR86d:47057} that the radical of the Volterra
subalgebra of the Cuntz algebra is the closed commutator ideal of the
Volterra subalgebra.  This result extends to graph \cstar algebras.
In the proposition below, we let $\sC$ denote the closed ideal
generated by the commutators of $\Alg \sN$.  The proof differs from
the one in~\cite{MR86d:47057}, which does not use groupoid techniques.

\begin{proposition} \label{p:radical}
The radical of $\Alg \sN$ is equal to the closed commutator
ideal $\sC$.
\end{proposition}

\begin{proof}
As usual, we  view all elements of $\sA$ as
functions on $\sG$.  The multiplication in $\sA$ is then given by a
convolution formula.
By Proposition~\ref{spec_rad}, 
$R = \{(x,k,y) \in \sigma(\Alg \sN) \mid x \prec y \}$ is the spectrum
of the radical and $A(R)$ is the radical of $\Alg \sN$.

Let $f,g \in \Alg \sN$.  
If we show that $[f,g](x,k,x)=0$ whenever 
$(x,k,x) \in \sigma(\Alg \sN)$ then $[f,g] \in A(R)$ and
$\sC \subseteq A(R)$.
Now
\begin{displaymath}
 f\cdot g(x,k,x) = \sum f(x,i,u)g(u,k-i,x) 
\end{displaymath}
where the sum is taken over all $i \in \bbZ$ and $u \in P$ for
which $(x,i,u)$ and $(u,k-i,x)$ lie in $\sigma(\Alg \sN)$.  This
requires both $x \preceq u$ and $u \preceq x$, so the only possibility
for $u$ is $u=x$.  If we make the change of variable $j = k-i$, then 
\begin{align*}
f\cdot g(x,k,x) &= \sum_i f(x,i,x)g(x,k-i,x) \\
 &=\sum_j f(x,k-j,x)g(x,j,x) \\
&= g\cdot f(x,k,x)
\end{align*}
Thus, $[f,g]=fg-gf$ vanishes at all points in $\sigma(\Alg \sN)$ of
the form $(x,k,x)$ and so is supported on $R$.  This shows that
$[f,g] \in A(R)$ and it follows immediately that
$\sC \subseteq A(R)$.

To show that $A(R) \subseteq \sC$, it suffices, by the spectral
theorem for bimodules, to show that each Cuntz-Krieger partial
isometry from $A(R)$ is in $\sC$.  If 
$\SSadj{\alpha}{\beta} \in A(R)$, then there is $j$ such that
$\alpha_i = \beta_i$ for $i < j$ and $\alpha_j \prec \beta_j$.
The range projection for $\SSadj{\alpha}{\beta}$ is a
subprojection of $R_{\alpha_1 \dots \alpha_j}$ and the initial
projection is a subprojection of
$R_{\beta_1 \dots \beta_j}$.  Let $P$ be the smallest projection
in $\sN$ which contains $R_{\alpha_1 \dots \alpha_j}$.  Since
$R_{\alpha_1 \dots \alpha_j} \ll R_{\beta_1 \dots \beta_j}$, 
$R_{\beta_1 \dots \beta_j} \perp P$ and
$\SSadj{\alpha}{\beta} = P \SSadj{\alpha}{\beta} P^{\perp}$.  This
implies that 
\begin{displaymath}
  \SSadj{\alpha}{\beta} = P\SSadj{\alpha}{\beta} -
  \SSadj{\alpha}{\beta}P
= [P, \SSadj{\alpha}{\beta}] \in \sC
\end{displaymath}
and the proposition is proven.
\end{proof}

If we let $D = \{(x,0,x)\mid x \in P\}$, then 
$\sD =A(D) = \Alg \sN \cap (\Alg \sN)^*$.  
Since $D \cup R$ is a proper
subset of $\sigma(\Alg \sN)$, it follows that the norm closure of
$A(D) + A(R)$ is a proper subset of $\Alg \sN$.  Thus, 
$\Alg \sN$ does not have a radical plus diagonal decomposition.
Furthermore, since $\sigma(\Alg \sN) \cup \sigma(\Alg \sN)^{-1}$ is a
proper subset of $\sG$, $\Alg \sN + (\Alg \sN)^*$ is not norm dense in
$\sA$.  This says that $\Alg \sN$ is ``non-Dirichlet.''  When every
loop has an entrance, $\sD$ is a masa in $\sA$ and $\Alg \sN$ is
triangular, but not strongly maximal triangular.
 However, we do have the following Proposition.

\begin{proposition} \label{p:maxtriang}
Assume that $G$ is a finite graph in which every loop has an entrance.
$\Alg \sN$ is maximal triangular in $\sA$.
\end{proposition}

\begin{proof}
Since $\Alg \sN \cap (\Alg \sN)^* = \sD$ is a masa, $\Alg \sN$ is a
triangular subalgebra of $\sA$.
Assume that $\Alg \sN \subseteq \sT \subset \sA$ and that $\sT$
is triangular.  It follows that $\sT \cap {\sT}^* = \sD$.  Let
$T \in \sT$.  Let $P$ be a projection in  $\sN$.
 Since $PTP^{\perp}$ leaves invariant each
projection in $\sN$, $PTP^{\perp} \in \Alg \sN$. 
Since $\sN \subseteq \sT$, $P^{\perp}TP \in \sT$.  It follows that
$S = PTP^{\perp} + P^{\perp} TP$ is a self-adjoint element of $\sT$
and hence lies in $\sD$.  Since $\sD$ is abelian and $P \in \sD$, 
$P$ commutes with $S$.  Thus $0 = P^{\perp} SP$.  But
$P^{\perp} SP = P^{\perp} TP$, so $P^{\perp} T P =0$.  Thus $T$ leaves
invariant each projection in $\sN$ and so must be an element of
$\Alg \sN$.  This shows that $\sT \subseteq \Alg \sN$ and $\Alg \sN$
is maximal triangular.
\end{proof}

\section{Normalizing Partial Isometries}\label{s:normalizing}

In this section we
characterize the partial isometries in a graph C*-algebra 
$\sA = C^*(G)$ which normalize the canonical diagonal algebra $\sD$.
We  assume throughout that  $G$ is 
a countable range finite directed graph with no sources
such  that each loop
 has an entrance.
In particular, by  Theorem~\ref{t:masa}, this ensures that $\sD$ is a masa.
The characterization of the $\sD$-normalizing partial isometries will
be applied in section~\ref{s:triangular} to show that gauge invariant
triangular subalgebras are classified by their spectra.

Recall that a partial isometry $v$ is $\sD$-normalizing
if $v^*\sD v \subseteq \sD$ and $v\sD v^* \subseteq \sD$. We write $N_{\sD}(\sB)$
for the set of all $\sD$-normalizing partial isometries in a subset 
$\sB$. 
Also we say that $v_1 + \dots + v_n$ is an orthogonal sum of partial
isometries if the set of initial projections $v_i^* v_i$, and also
the set of final projections $v_iv_i^*$, consists of pairwise
orthogonal projections.

In Theorem~\ref{t:normalizing} we show that $\sD$-normalizing 
partial isometries are,
modulo coefficients from $\sD$, orthogonal sums of Cuntz-Krieger
partial isometries; moreover they are characterized by a property
which is preserved by isometric isomorphism. 
The equivalence of \ref{pinorm} and \ref{orthsum} in the case of
Cuntz algebras was obtained in \cite[Lemma~5.4]{MR99g:47108}, where it 
formed the basis for the calculation of normalizing partial isometry
homology groups of various triangular subalgebras.

\begin{theorem}\label{t:normalizing}
Let $v$ be a partial isometry in $\sA$.
Then the following assertions are equivalent\/{\rm :}
\begin{enumerate}
\item \label{pinorm} $v$ is a $\sD$-normalizing partial isometry. 
\item \label{orthsum} $v$ is an orthogonal sum of a finite number of partial 
isometries of the form $d \SSadj{\alpha}{\beta}$, where $d \in \sD$.
\item \label{zero_one} For all projections $p,q \in \sD$, the norm $\|qvp\|$ is 
equal to $0$ or $1$.
\end{enumerate}
\end{theorem}

This theorem is in complete analogy with
the following counterpart for 
AF C*-algebras. (See~\cite{MR91e:46078} 
or~\cite[Lemma~5.5 and Proposition~7.1]{MR94g:46001}.)
The Cuntz-Krieger partial isometries play the same role for graph 
\cstar algebras as systems of matrix units do for AF \cstar algebras.
 
\begin{theorem}\label{t:normalizingAF}
Let $\sB$ be an  AF \cstar algebra with finite dimensional subalgebra chain
$\sB_1 \subseteq \sB_2 \subseteq \dots $ and masas $\sC_k \subseteq \sB_k$
such that $N_{\sC_k}(\sB_k) \subseteq N_{\sC_{k+1}}(\sB_{k+1})$, for
all $k$. Suppose also that the union of the $\sB_k$ is dense
in $\sB$.
Then the closed union $\sC$
of the masas $\sC_k$ is a masa in $\sB$ and the following assertions are 
equivalent
for a partial isometry $v$  in $\sB$.
\begin{enumerate}
\item \label{afpinorm} $v$ is a $\sC$-normalizing partial isometry.
\item \label{aforthsum} $v = cu$ where $c \in \sC$ and 
$u \in  N_{\sC_k}(\sB_k)$, for some $k$.
\item \label{afzero_one} For all projections $p,q \in \sC$, 
the norm $\|qvp\|$ is equal to $0$ or $1$.
\end{enumerate}
\end{theorem}

Theorem~\ref{t:normalizingAF} will be used in the proof of 
Theorem~\ref{t:normalizing}
to show that if $v$ is in $N_{\sD}(\sA)$ then so too is its AF part
$v_0 = \Phi_0(v)$.
We  also require the following two lemmas.

\begin{lemma} \label{l:path_proj}
Let $\alpha$, $\beta$  be paths of the same length and
let $e=\SSadj{\alpha}{\beta}$ be a non-zero partial isometry in
the AF subalgebra $\sF$ of $C^*(G)$.  
For each positive integer $k$  there exist non-zero
projections $q$, $p$ with $q = e p e^*$, such that for all 
paths $\gamma$ with length at most $k$, and for all
$\SSadj{\lambda}{\mu}$ in $\sF$ with 
$| \lambda |= | \mu |\leq k$, we have 
\[
q(S_{\gamma}^{\vphantom{*}}\SSadj{\lambda}{\mu})p =
q(\SSadj{\lambda}{\mu}S_{\gamma}^{\vphantom{*}})p = 0.
\]
\end{lemma} 

\begin{proof}
Note that if we verify the lemma for an integer $k$, then we have
verified it for all integers less than $k$; thus we may increase a
value for $k$ if needed.  This, together with the hypothesis that
every loop has an entrance, allows us to choose a path
$\pi = f_{2k} \dots f_1$ of length $2k$ and a path 
$w = w_1 w_2 \dots $ of length at least $k$ such that
\begin{enumerate}
\item $r(\pi) = r(f_{2k}) = s(\alpha) = s(\beta)$.
\item $r(w) = s(\pi) = s(f_1)$.
\item For every integer $d$ with $1 \leq d \leq k$, 
$f_d \dots f_1 \neq w_1 \dots w_d$.
\end{enumerate}
Indeed, the assumption that there are no sources allows us to choose
the path $\pi$ with range vertex equal to $s(\alpha)=s(\beta)$.
Possibly, we can choose $\pi$ so that
$r(f_d) \neq s(f_1)$ for $d = 1, \dots, k$.
In this case, any extension $w = w_1 \dots w_k$ works, since
$r(w_1 \dots w_d) = r(w_1)= s(f_1) \neq r(f_d) = r(f_d \dots f_1)$.  
On the other hand, if we must back into a loop then (increasing $k$ if
necessary), we can arrange that $f_k \dots f_1$ consists of multiple
repeats of a single simple loop.  By choosing $w$ so that $w_1$ is an
entrance to the loop, we guarantee that
$f_d \dots f_1 \neq w_1 \dots w_d$ for all $d \leq k$.

Now let
\[
p = \ssadj {\beta \pi w} \text{ and } q = \ssadj {\alpha \pi w} = epe^*.
\]
Let $\gamma = \gamma_1 \dots \gamma_d$ be a path with 
$1 \leq d = |\gamma| \leq k$.  We first show that
$q S_{\gamma} p =0$.  If not, then
\[
\ssadj {\alpha \pi w} S_{\gamma}^{\vphantom{*}} \ssadj {\beta \pi w}
\neq 0.
\]
Now, recall that for any edges $e$ and $f$, $\SadjS{e}{f} = 0$ except
when $e=f$ (the ranges of the generating partial isometries $S_e$ are
pairwise orthogonal) and that, if $g$ is an edge with
$r(g) = s(e)$ then $\sadjs{e}S_g^{\vphantom{*}} = S_g$ (from the
Cuntz-Krieger relations).  Consequently, the edges in the finite
path $\alpha \pi w$ match the edges (reading from left to right) in
the path $\gamma \beta \pi w$.  Since the length of $\pi$ is at least
twice the length of $\gamma$, the cancellations into $\gamma \beta \pi
w$ bring us $d$ edges into $w$; this forces
$f_d \dots f_1 = w_1 \dots w_d$.  But this contradicts the choice of
$\pi$ and $w$.

Insertion of $\SSadj {\lambda}{\mu}$ either before or after
$S_{\gamma}$ does not change the result: there are at most $k$
cancellations from $S_{\mu}^*$, which cannot affect the second half of
$S_{\pi}$ since $|\pi| = 2k$, and the cancelled partial isometries are
replaced by partial isometries from $S_{\lambda}$.  Thus the general
result holds.
\end{proof}

In Lemma~\ref{l:proj_afpart}, $B^*(G)$ is the (non-closed) algebra
generated by the Cuntz-Kreiger partial isometries and the maps
$\Phi_m$ are as defined in section~\ref{s:prelim}.

\begin{lemma} \label{l:proj_afpart}
Let $a \in B^*(G)$ and let $e = \SSadj{\alpha}{\beta}$ be a partial
isometry in the AF subalgebra $\sF$ of $C^*(G)$. 
 Then there exist projections $p$
and $q = epe^*$ such that 
$qap = q\Phi_0(a) p$.
\end{lemma}

\begin{proof}
By the observations in section~\ref{s:prelim}, $a - \Phi_0(a)$ can be
written as a finite linear combination of terms of one of the two
forms: $S_{\gamma}^{\vphantom{*}} \SSadj {\lambda}{\mu}$ and
$\SSadj {\lambda}{\mu} S_{\gamma}^*$, where $|\gamma| \geq 1$ and
$|\lambda| = |\mu|$.  An application of Lemma~\ref{l:path_proj} gives
projections $p_1$ and $q_1$ such that $q_1 = e p_1 e^*$ and
$q_1 S_{\gamma}^{\vphantom{*}} \SSadj {\lambda}{\mu} p_1 = 0$ 
for all terms of this type in the linear combination.  Now let 
$f = (ep_1)^*$ and apply Lemma~\ref{l:path_proj} again to obtain
projections $q \leq q_1$ and $p \leq p_1$ with
$f q f* = p$ and 
$p S_{\gamma}^{\vphantom{*}} \SSadj {\mu}{\lambda} q=0$ for all terms
of the second type in the linear combination for $a - \Phi_0(a)$.  It
now follows that 
$q \SSadj {\lambda}{\mu} S_{\gamma}^* p = 0$ and
$q S_{\gamma}^{\vphantom{*}} \SSadj {\lambda}{\mu} p = 0$ for all
the terms; hence $qap = q \Phi_0(a) p$.
\end{proof}

\begin{proof}[Proof of Theorem~\ref{t:normalizing}]
The implications~(\ref{orthsum})~$\Longrightarrow$~(\ref{pinorm}) 
and~(\ref{pinorm})~$\Longrightarrow$~(\ref{zero_one}) are routine.  

It remains prove
that~(\ref{zero_one})~$\Longrightarrow$~(\ref{orthsum}):
Let $v$ be a partial isometry in $\sA$ which satisfies
condition~(\ref{zero_one}).  We claim first that
$\Phi_0(v)$ is a $\sD$-normalizing partial isometry. If not,
then, since $\Phi_0(v)$ is in the AF subalgebra $\sF$, we can use 
Theorem~\ref{t:normalizingAF} to find a partial isometry
$e=\SSadj{\alpha}{\beta}$ in $\sF$ and a $\delta >0$ such that, for
any pair $p \leq e^*e$, $q = epe^*$ we have
\[
\delta \leq \|q\Phi_0(v)p\| \leq 1-\delta.
\]
(This is easy to do using the function representation of $v_0$ on the
AF subgroupoid and knowledge of the form of $v_0$ -- that it is not an
element of $\sD$ times a matrix unit.)

Now let $v' \in B^*(G)$ be such that
$\|v'-v\| < \delta/2$.  By Lemma~\ref{l:proj_afpart}, there exist
non-zero projections $p$ and $q$ with $q=epe^*$ such that
$q\Phi_0(v')p = qv'p$. Since $\|qvp\|$ is either 0 or 1, either
$\|qv'p\| \leq \delta/2$ or $1-\delta/2 \leq \|qv'p\|$.
Since $\|\Phi_0(v') - \Phi_0(v)\| < \delta/2$ ($\Phi_0$ is
contractive), this implies that either
$\|q \Phi_0(v) p \| < \delta$ or $1-\delta < \|q\Phi_0(v)p\|$, a
contradiction.   Thus, the 0-order term in the `Fourier' series for
$v$ is $\sD$-normalizing.  It follows from Theorem~\ref{t:normalizingAF}
that $\Phi_0(v)$ has the form required in condition~(\ref{orthsum}).

Now suppose that $m>0$.  If $|\nu|=m$ and
$|\lambda|-|\mu|=m$, then the product
$S_{\nu}^*\SSadj{\lambda}{\mu}$ is either zero or of the form
$\SSadj{\lambda_1}{\mu}$ with $|\lambda_1|=|\mu|$.  
It follows that 
$S_{\nu}^* \Phi_m(v) = \Phi_0(S_{\nu}^*v)$.  Since $v$ satisfies
condition~(\ref{zero_one}), so does $S_{\nu}^*v$; the argument above
shows that $S_{\nu}^* \Phi_m(v)$ is $\sD$-normalizing and has the required form
(condition~(\ref{orthsum})).  This, in turn, implies that
$\ssadj{\nu}\Phi_m(v)$ is $\sD$-normalizing and has the required form
for every path $\nu$ with length $m$.  Consequently, $\Phi_m(v)$
satisfies condition~(\ref{orthsum}).  In a similar fashion, we can
show that when $m<0$, $\Phi_m(v)$ satisfies condition~(\ref{orthsum})
(consider adjoints, for example).

Now, if $w$ is a partial isometry and $ww^*xw^*w \neq 0$ then
$\|w + ww^*xw^*w\| > 1$.  From this observation and the Cesaro
convergence, it follows that the operators $\Phi_m(v)$ are non-zero
for only finitely many values of $m$ and that $v$ is the orthogonal
sum of these operators.  Thus $v$ has the form required in 
condition~(\ref{orthsum}).
\end{proof}

\section{Triangular Subalgebras Determine Their Spectrum}
\label{s:triangular}

In this section we show that the gauge invariant triangular subalgebras
of certain graph C*-algebras are classified by the
isomorphism type of their spectra.
We  assume throughout the section that 
$G_1$ and $G_2$ are 
countable range finite directed graphs with no sources
and   that each loop
 has an entrance.

\begin{theorem} \label{tri_class}
For $i =1,2$, let $\sT_i$ be a triangular subalgebra of $\sA_i$ with
diagonal $\sD_i$ such that $\sT_i$ generates $\sA_i$ as a
 \cstar algebra and $\sT_i$ is invariant under the gauge
 automorphisms \/{\rm (}so that $\sT_i = A(\sP_i)$, where
$\sP_i = \sigma(\sT_i)$\/{\rm )}.
Then the following statements are equivalent\/{\rm :}
\begin{enumerate}
\item \label{algiso} $\sT_1$ and   
 $\sT_2$ are isometrically
isomorphic operator algebras.
\item \label{speciso}  There is a groupoid isomorphism 
$\gamma \colon \sG_1 \to  \sG_2$
with    $\gamma (\sP_1) = \sP_2$.
\end{enumerate}
\end{theorem}

\begin{proof}
If $\gamma$ has the properties of (\ref{speciso}) it is plain
that there is a \cstar algebra isomorphism
$C^*(\sG_1) \to C^*(\sG_2)$
which restricts to an isometric
isomorphism
$\sT_1 \to \sT_2$.
Assume then that  $\Gamma \colon \sT_1 \to \sT_2$
is an isometric isomorphism.
In view of the norm characterization of normalizing partial isometries,
Theorem~\ref{t:normalizing} condition~(\ref{zero_one}), we have
$\Gamma (N_{\sD_1}(\sT_1)) =  N_{\sD_2}(\sT_2)$.
Moreover, since $\Gamma(\sD_1) = \sD_2$,
the groupoid support of $\Gamma (du)$ with $d \in \sD_1$
and $u$ in $N_{\sD_1}(\sT_1)$ is independent
of $d$ if $d$ is a partial
isometry and $d^*du = u$.
Reciprocally, in view of Theorem~\ref{t:normalizing} condition~(\ref{orthsum}),
normalizing partial isometries are determined by their
groupoid support, up to a diagonal multiplier.
Thus $\Gamma$ induces a
map
\[
\tilde{\gamma} : \{\supp(u): u \in
N_{\sD_1}(\sT_1)\} \to \{\supp(v): v \in N_{\sD_2}(\sT_2)\}.
\]

Each point $g$ in $\sP_1$ is an intersection of the
sets $Z(\alpha,\beta)$ that contain it and so by
local compactness  $\tilde{\gamma}$ defines
a bijection $\gamma \colon \sG_1 \to \sG_2$
and this map, in turn, induces   $\tilde{\gamma}$.
In particular, $\gamma$ is a homeomorphism.

Note now that  $\gamma$ is a semigroupoid map. 
To see this observe first that if $u_1$ and $u_2$
are normalizing partial isometries in $\sT_1$ with 
support sets $U_1$ and $U_2$ then $u_1u_2$ has support set $U_1 \cdot U_2$.
Thus, if $g_1$ and $g_2$ are composable elements in  $\sP_1$ and
 \begin{displaymath}
\{g_1\} = \cap_{n=1}^\infty U_n,\quad  \{g_2\} = \cap_{n=1}^\infty V_n,
 \end{displaymath}
where $U_n, V_n$ are supports of normalizing partial isometries,
then
 \begin{displaymath}
\{g_1 \cdot g_2\} = \cap_{n=1}^\infty U_n\cdot V_n.
 \end{displaymath}
Thus
 \begin{displaymath}
 \{ \gamma(g_1 \cdot g_2)\} = 
\cap_{n=1}^\infty \tilde{\gamma}(U_n\cdot V_n)
 = \cap_{n=1}^\infty \tilde{\gamma}(U_n)\cdot \tilde{\gamma}(V_n)
 \end{displaymath}
and this last set is the singleton set
$\{\gamma(g_1) \cdot \gamma(g_2)\}$.

We  now extend $\gamma$ to a map from $\sG_1$
to $\sG_2$.
Note first  that since, by hypothesis,
$\sT_1$ generates $\sA_1$ as a C*-algebra,
the sets
\[
U = U_1 \cdot U_2^{-1} \cdot U_3 \cdot U_4^{-1} \cdot \ldots 
\cdot U_{2n}^{-1},
\]
where $U = \supp(u_i)$ and
$u_i \in N_{\sD_1}(\sT_1)$,
have union equal to
$\sG_1$.
Indeed, approximate $u$ in $N_{\sD_1}(\sA_1)$ by a polynomial
$p$ in the generators and their adjoints,
\begin{displaymath}
  p = \sum_m \sum_{|\lambda|-|\mu|=m} a_{\lambda \mu}
\SSadj{\lambda}{\mu}.
\end{displaymath}
and it follows that
a point $g$ in $\supp(u)$ must lie
in the support of some $\SSadj{\lambda}{\mu}$
in the sum.
Extend $\gamma$ to $\sG_1$ by setting
\[
\gamma(g_1 \cdot g_2^{-1} \cdot \ldots \cdot g_{2n}^{-1})
=  \gamma(g_1) \cdot \gamma(g_2^{-1}) \cdot \ldots \cdot
      \gamma(g_{2n}^{-1}).
\]
This is well defined and onto, since $\sT_2$ generates 
$\sA_2$, and so, as before, this extension   is  a groupoid
isomorphism.
\end{proof}

It is clear that the proof method above simplifies to give
the following equivalence between  the isomorphism type 
of  the pair
$(C^*(G),\sD)$ and the isomorphism type of
the  groupoid $\sG$.
(Compare, for example, \cite[Theorem~7.5]{MR94g:46001}.)

\begin{theorem}  
The following statements are equivalent\/{\rm :}
\begin{enumerate}
\item  There is a C*-algebra isomorphism
$\Gamma : C^*(G_1) \to C^*(G_2)$ with
$\Gamma(\sD_1) = \sD_2$,
where each $\sD_i$ is the canonical abelian diagonal
subalgebra of $C^*(G_i)$.
\item There is a groupoid isomorphism
$\gamma \colon \sG_1 \to \sG_2.$
\end{enumerate}
\end{theorem}

Similarly the proof above extends with only trivial changes
to give an equivalence between the isomorphism
type of the pair $(\sB,\sD)$, consisting of a gauge invariant
subalgebra $\sB$ containing the diagonal $\sD$, and the isomorphism type of 
the spectrum of $\sB$.


\providecommand{\bysame}{\leavevmode\hbox to3em{\hrulefill}\thinspace}
\providecommand{\MR}{\relax\ifhmode\unskip\space\fi MR }
\providecommand{\MRhref}[2]{%
  \href{http://www.ams.org/mathscinet-getitem?mr=#1}{#2}
}
\providecommand{\href}[2]{#2}

\end{document}